\documentclass[preprint,3p]{elsarticle}
\makeatletter 
\@mparswitchfalse%
\makeatother
\normalmarginpar

\usepackage{todonotes}
\makeatletter
\renewcommand{\todo}[2][]{\tikzexternaldisable\@todo[#1]{#2}\tikzexternalenable}
\makeatother

\usepackage[verbose]{placeins}
\usepackage{bmk_maths}

\usepackage{tikz}
\usepackage{graphicx} 
\usepackage{subcaption}
\usepackage{import}
\usepackage{currfile}
\biboptions{numbers,sort&compress} 

\usepackage{pgfplots}
\pgfplotsset{compat=newest}
\usepackage{pgfplotstable}
\usepgfplotslibrary{polar}
\usepgfplotslibrary{statistics}
\usetikzlibrary{external}
\usetikzlibrary{spy}
\usepgfplotslibrary{patchplots}
\usepackage{bmk_plots}

\usepackage{xstring}
\renewcommand{\inputfrom}[2]{%
  \StrSubstitute{#2}{.tex}{}[\cleanfilename]%
  \StrSubstitute{#1\cleanfilename}{/}{-}[\temp]%
  \includegraphics{\temp.pdf}%
}

\usepackage{amsmath}
\usepackage{amsthm}
\usepackage{amssymb}
\usepackage{amsfonts}
\usepackage{matlab-prettifier}

\usepackage{pgfplotstable}
\usepackage{booktabs}
\usepackage{algorithm}
\usepackage{algpseudocode}  
\usepackage{hyperref}
\usepackage{cleveref}
\Crefname{figure}{Figure}{Figures}

\newdefinition{remark}{Remark}

\Crefname{assumption}{Assumption}{Assumptions}
\newtheorem{example}{Example}
\Crefname{example}{Example}{Examples}

\makeatletter
\newcommand{\algmargin}{\the\ALG@thistlm}
\makeatother
\newlength{\whilewidth}
\algdef{SE}[parWHILE]{parWhile}{EndparWhile}[1]
  {\parbox[t]{\dimexpr\linewidth-\algmargin}{%
     \hangindent\whilewidth\strut\algorithmicwhile\ #1\ \algorithmicdo\strut}}{\algorithmicend\ \algorithmicwhile}%
\algnewcommand{\parState}[1]{\State%
  \parbox[t]{\dimexpr\linewidth-\algmargin}{\strut #1\strut}}

\newcommand{\AH}[1]{{#1}}
\newcommand{\ah}[1]{#1}
\newcommand{\bk}[1]{#1}
\newcommand{\lt}[1]{#1}
\newcommand{\mg}[1]{#1}
\newcommand{\MGnew}[1]{{#1}}

\usepackage[normalem]{ulem} 
\newcommand{\bkrevise}[1]{{#1}}
\newcommand{\bkrevisesout}[2]{{#2}}
\newcommand{\mgrevise}[1]{{#1}}

\journal{Journal of Computational Physics}

\begin{document}

\begin{frontmatter}



    \title{Noise-robust multi-fidelity surrogate modelling for parametric partial differential equations}


    \author[1]{Benjamin M. Kent\corref{cor1}}
    \ead{kent@imati.cnr.it}
    \author[1]{Lorenzo Tamellini}
    \author[LaCaN,CIMNE]{Matteo Giacomini}
    \author[LaCaN,CIMNE]{Antonio Huerta}
    \cortext[cor1]{Corresponding author}

    \affiliation[1]{organization={Consiglio Nazionale delle Ricerche -- Istituto di Matematica Applicata e Tecnologie Informatiche ``E. Magenes'' (CNR--IMATI)},
        addressline={Via Adolfo Ferrata, 5},
        city={Pavia},
        postcode={27100},
        state={PV},
        country={Italy}}
    \address[LaCaN]{Laboratori de C\`alcul Num\`eric (LaC\`aN), ETS de Ingenier\'ia de Caminos, Canales y Puertos, Universitat Polit\`ecnica de Catalunya - BarcelonaTech (UPC), Barcelona, Spain.}
    \address[CIMNE]{Centre Internacional de M\`etodes Num\`erics en Enginyeria (CIMNE), Barcelona, Spain.}

    \begin{abstract}
        We address the challenge of constructing noise-robust surrogate models for quantities of interest (QoIs) arising from parametric partial differential equations (PDEs), using multi-fidelity collocation techniques; specifically, the Multi-Index Stochastic Collocation (MISC).
        In practical scenarios, the PDE evaluations used to build a response surface are often corrupted by numerical
        noise, especially for the low-fidelity models.
        This noise, which may originate from loose solver tolerances, coarse discretisations, or transient effects, can lead to overfitting in MISC, degrading surrogate quality through nonphysical oscillations and loss of convergence, thereby limiting its utility in downstream tasks like uncertainty quantification, optimisation, and control.
        To correct this behaviour, we propose an improved version of MISC that can automatically detect the presence of solver noise during the surrogate model construction and then ignore the exhausted fidelities.
        Our approach monitors the spectral decay of the surrogate at each iteration, identifying stagnation in the coefficient spectrum that signals the onset of noise. Once detected, the algorithm selectively halts the use of noisy fidelities, focusing computational resources on those fidelities that still provide meaningful information.
        The effectiveness of this approach is numerically validated on two challenging test cases: a parabolic
        advection--diffusion PDE with uncertain coefficients, and a parametric turbulent incompressible Navier--Stokes problem. The results showcase the accuracy and robustness of the resulting multi-fidelity surrogate and its capability to extract relevant information, even from under-resolved meshes not suitable for reliable single-fidelity computations.
    \end{abstract}




\end{frontmatter}





\section{Introduction}
Parametric partial differential equations (PDEs) are commonly used to model physical problems with \mg{user-defined, possibly uncertain, configurations}.
The dependence of the solution on such \mg{(uncertain)} parameters, or that of a quantity of interest derived from the solution, is often well approximated using global polynomial approximation under mild assumptions on the parameter dependence \lt{of} the problem.
Popular constructions include \bk{using} sparse grids, radial basis functions, Gaussian processes, least squares methods, stochastic Galerkin, reduced basis methods, \mg{proper orthogonal and} proper generalized decompositions and neural networks \cite{ghanem:uqbook,Smith2013,Sullivan2015,gramacy2020surrogates,rozza:book,quarteroni:RBbook,Adcock2022,kutiniok:DNN,hesthaven:nn,fresca.eal:nn,adcock.dexter:DNN,Chinesta2017,Giacomini2021}.
Many of these constructions require repeatedly solving the PDE.
\bk{To reduce computational costs, so-called multi-fidelity variants of these methods have been introduced.
    These may use greedy adaptive algorithms in both parametric and physical \ah{(space, time...)} domains \cite{Chkifa2014,Eigel2014,Nobile2016,Bespalov2016}, or multi-mesh and multi-level approaches \cite{Jakeman2019,HajiAli2016b,Crowder2019,Bespalov2022b,legratiet:multifidgp,Ng2012};
    see e.g., \cite{peherstorfer:mfsurvey} for a survey.}

\lt{This work concerns multi-fidelity approximation. Besides the already-hinted mesh or timestep refinement, a multi-fidelity approach could also consider solvers with
    different tolerances, or even a set of different physical models of the system at hand.}
Specifically, we focus on \emph{multi-index stochastic collocation} \mg{(MISC)} as introduced in \cite{HajiAli2016a}. This is a multi-fidelity surrogate modelling technique that builds upon the well-known sparse grid stochastic collocation (\lt{SGSC}) approximation method for parametric PDEs \cite{Babuska2010,Xiu2005,Gerstner2003,Piazzola2022}.
Methodologies built upon \lt{SG}SC are popular due to their non-intrusive nature: a quantity of interest \mg{needs to be evaluated only} for a prescribed collocation set of parameter realisations, and the surrogate model is then built by interpolating these evaluations.
In \lt{our multi-fidelity} \mg{context}, the collocation set specifies both \lt{parameter realisations} and solver fidelities, \mg{which} are combined in a similar way to \lt{SG}SC methods to create a multi-fidelity surrogate model.
Adaptive MISC techniques have been successfully applied in a number of practical applications, see e.g., \cite{Jakeman2019,Piazzola2022a,Chiappetta2024}.

Adaptive \lt{SG}SC and MISC algorithms generally use a ``greedy'' approach that refines the surrogate model by seeking the greatest change to the computed surrogate.
This assumes that every change in the surrogate is a feature that must be captured to \mg{achieve} an accurate final approximation.
In many applications this is \lt{however} not the case: the underlying solvers (especially the lowest fidelities) may introduce spurious, uncontrolled errors \lt{(e.g., \mg{due to} iterative solver error, discretisation errors,...)}.
\mg{Hence} greedy \mg{strategies} may \mg{be unable to} distinguish between the \emph{true} solution behaviour and the solver error, \mg{and consequently ``overfit'' to the inherent solver noise}.
This \mg{leads to} spurious oscillations in the resulting surrogate models.
See for example \cite{Piazzola2022a,brown.back:noise,deBaar2015,pellegrini.eal:noisySRBF}, as well as \cite{lukaczyk:noise,geraci:acv.nondet} and references therein for
further discussion on noise corruption in the context of outer loop applications such as \bk{uncertainty quantification} (UQ) or shape optimisation.

\mgrevise{
Classical least-squares approaches (either single- or multi-fidelity) to surrogate modelling construction can mitigate the impact of noisy evaluations of the quantity of interest \cite{wolfers:multi-level-PCE,migliorati:rdpforfunc,migliorati:rdpforpde,sudret:adaptive,hampton.doostan:coherence,narayan.eal:christoffel}, although this noise is not always explicitly accounted for in these approaches; the impact of noise in the convergence of least squares approximations is discussed in e.g.\ \cite{coehn.migliorati:noisy,cohen.leviatan:noisyLS,kohler:noisyLS,gyorfi:noisyLS}.
Among the alternative approaches for dealing with noise in model evaluations, we mention stochastic radial basis functions \cite{pellegrini.eal:noisySRBF,Piazzola2022a}, the approximate control variate method for non-deterministic models discussed in \cite{geraci:acv.nondet}, and the noise-robust approximation of level sets of multivariate functions discussesd in \cite{croci2025:noisylevelset}.

With the same purpose, in this work we propose instead a novel noise-robust variant of the MISC algorithm, that we call \plateaumisc{}. 
\plateaumisc{} modifies the greedy MISC algorithm to differentiate between refinements that capture the \emph{true} parametric response and refinements that are capturing solver noise. 
\lt{W}e achieve this by inspecting the spectral content of the approximation \mg{and automatically detecting deviations of the behavior of the multi-fidelity response surfaces from the expected spectral decay of smooth quantities of interest}, see e.g.\ \cite[Chapter 4]{Adcock2022}.
\mg{More precisely, the algorithm \mg{detects} the corruption due to solver error by a failure to adhere to an a priori bound for the spectral polynomial expansion coefficients.}

Compared to the classical version of MISC, the \plateaumisc{} algorithm thus offers a multi-fidelity approximation with the following advantages:
\begin{itemize}
    \item it is robust to solver noise, allowing the exploitation of information from low fidelity approximations without corrupting the final approximation;
    \item it automatically and efficiently transitions from low-fidelity to high-fidelity solvers, blocking the ``exhausted'' low fidelity solvers, as the computational budget increases.
\end{itemize}
}

\AH{
We conclude this introduction by briefly mentioning a related line of research on fault-tolerant surrogate model construction. This approach addresses settings in which a \emph{small} number of evaluations of the quantity of interest are severely \emph{corrupted} or even \emph{unavailable} due to simulation failures. This contrasts with the scenario considered in the present work, where we assume that a \emph{large} number of evaluations are affected by \emph{noise}. Fault-tolerant surrogate modeling techniques based on compressed sensing are discussed in \cite{adcock.jakeman:noise-fault,shin.xiu:fault-tolerant} and could, in principle, be adapted to the setting studied here, whereas the fault-tolerant sparse grid approximations proposed in \cite{heene:fault-tol-SG,harding:fault-tol-SG} do not extend to our scenario; see also \cite{pauli:fault-tolerant} for a fault-tolerant multi-level Monte Carlo method.
}

The paper is structured as follows.
\Cref{sec:surrogates} details the construction of a single-fidelity \lt{SG}SC approximation and the polynomial approximation properties that will be exploited.
\Cref{sec:pmisc} presents the classic MISC algorithm and \mg{the} novel \plateaumisc{} adaptive algorithm exploiting spectral coefficient plateaus.
\Cref{sec:numerics} \mg{numerically} demonstrates the algorithm \mg{applied to} a parabolic advection--diffusion PDE \mg{with uncertain parameters} and a \mg{parametric} Reynolds Averaged Navier--Stokes flow \mg{with Spalart--Allmaras turbulence model}.
The paper concludes in \Cref{sec:conclusions}.

\section{Collocation based methods for parametric PDEs}\label{sec:surrogates}
In abstract terms, the purpose of this work is to accurately and efficiently approximate a multi-variate function $q:\Gamma\subset\R^{\nY} \to \R$ whose values are only known implicitly upon solving a PDE.
The inputs $\y\in\Gamma$ will be referred to as parameters \ah{and we consider a Cartesian} parameter domain $\Gamma:=\Gamma_{1}\times\Gamma_2\times \cdots \times \Gamma_{\nY}$.
For simplicity it is assumed $\Gamma_{\iY}=[0,1]$ for $\iY=1,2,...,\nY$ and \lt{that} the quantity of interest \mg{(QoI)} $q$ is a scalar-valued function.
\mg{For instance, in the framework of uncertainty quantification} the parameters represent images of underlying random variables \mg{(\Cref{sec:parabolic}), whereas in parametric flow problems they describe working conditions of the system (\Cref{example:2dwmh}).}

\mg{Assuming \lt{that} our quantity of interest is smooth with respect to the parameters $\y$, global polynomial approximation represents a common choice in the context of parametric PDEs and related problems, as already \lt{mentioned} in the introduction.} 
In particular, this assumption ensures that we can write the quantity of interest as a convergent polynomial series.
\ah{More explicitly, consider a basis of multi-variate polynomials $\{\chebBasis_{\vMiPoly}(\y)\}_{\vMiPoly\in\N_0^{\nY}}$ where  $\chebBasis_{\vMiPoly}(\y) = \prod_{\iY=1}^{\nY}\oneChebBasis_{\miPoly_{\iY}}(y_{\iY})$ and each term $\oneChebBasis_{\miPoly_{\iY}}(y_{\iY})$ is a polynomial of degree $p_{\iY}$ \mg{in the parameter $y_\iY$}.
For example, $\oneChebBasis_{\miPoly_{\iY}}(y_{\iY})$ could be the monomial $y^{\miPoly_{\iY}}$, or the Chebyshev polynomial of the first kind of degree ${\miPoly_{\iY}}$.
}
The quantity of interest $q$ can then be expressed as a convergent polynomial series,
\begin{equation}\label{eq:spectral}
    q(\y) \equiv \sum_{\vMiPoly \in \N_0^{\nY}} \FT{q}_{\vMiPoly}\chebBasis_{\vMiPoly}(\y).
\end{equation}

In this work we consider global polynomial approximations constructed via sparse grid stochastic collocation (\lt{SG}SC) type methods which we now briefly recall, \lt{see e.g.\ \cite{Babuska2007,Xiu2005,Piazzola2022}}.
\lt{In particular, we closely follow the setup of \cite{Piazzola2022}.}
\paragraph{\ah{One dimensional interpolation}}
First, denote by $\{\collocSet{\level}\}_{\level\in\N}$ a sequence of sets of one dimensional collocation points where $\collocSet{\level}:=\{\collocPt_{\iColloc}\}_{\iColloc=1}^{\nColloc{\level}}\subset[0,1]$.
\bk{The number of points $\nColloc{\level}$ at each level $\level\in\N$ is specified by the so-called \emph{level-to-knots} function.}
Common choices of the level-to-knots function $\nColloc{\level}$ are
\begin{equation}\label{eq:rules}
    \nColloc{\level}=1 \text{ for } \level=1, \text{ and } \nColloc{\level}=\begin{cases}
        \level         & \text{`linear' \mg{rule}}   \\
        2(\level-1)+1  & \text{`two-step' \mg{rule}} \\
        2^{\level-1}+1 & \text{`doubling' \mg{rule}} \\
    \end{cases} \forall \level > 1.
\end{equation}
A common choice of points $\collocSet{\level}:=\{\collocPt_{\iColloc}\}_{\iColloc=1}^{\nColloc{\level}}$ are the Clenshaw--Curtis (CC) points which are defined by the extrema of the Chebyshev polynomials \cite{Clenshaw1960}.
CC points are nested under the `doubling' level-to-knots rule of \eqref{eq:rules}.
In this work, we also use the symmetric Leja points for a uniform random variable.
By construction, Leja points are a nested sequence of points for every choice of \bk{level-to-knots rule} in \eqref{eq:rules} \cite{Piazzola2022}.
\ah{Clenshaw--Curtis points, symmetric Leja points and the `doubling' and `two-step' rules are illustrated in \Cref{fig:knots}.}

Next, we define a corresponding sequence of interpolation operators $\{\interpOp{\level}\}_{\level\in\N}$ such that for continuous $q\in C^0([0,1];\R)$ \mg{it holds that}
\begin{equation} \label{eq:oneDinterp}
    \interpOp{\level}[q](y):=\sum_{\collocPt\in\collocSet{\level}}q(\collocPt) \LagrangePoly{\collocPt}{\collocSet{\level}}(y),
\end{equation}
where $\{\LagrangePoly{\collocPt}{\collocSet{\level}}(y)\}_{\collocPt\in\collocSet{\level}}$ are the corresponding Lagrange interpolation polynomials $\LagrangePoly{\collocPt}{\collocSet{\level}}(y)=\prod_{w\in\collocSet{\level}\setminus\{\collocPt\}}(y-w)/(\collocPt-w)$.


\ah{\begin{figure}
        \inputfrom{figures/}{knots.tex}
        \caption{Example sequences of one-dimensional collocation points $\collocSet{\level}$ on $[0,1]$ for Clenshaw--Curtis points with the `doubling' rule and symmetric Leja points with the `two-step' rule, as defined in \eqref{eq:rules}.}
        \label{fig:knots}
    \end{figure}}

\paragraph{\ah{Multi-dimensional interpolation}}
A multi-dimensional approximation is constructed using \mg{the} multi-index set $\miset \subset \N^{\nY}$ to define a linear combination of tensor products of interpolation operators.
\lt{This is known to be an effective approach for approximating suitably smooth high-dimensional functions ($\nY \gg 1$) without resorting to tensor grid sampling of the parameter domain $\Gamma$ which quickly becomes computationally \mg{unaffordable} (i.e.\ \emph{the curse of dimensionality}) \cite{Barthelmann2000}.}
An approximation of a multi-variate function $q\in C^0(\Gamma;\R)$ is defined via the \emph{combination technique} formulation \cite{Griebel1992,Wasilkowski1995},
\begin{equation} \label{eq:combinationtechnique}
    q(\y) \approx \approxP{q}(\y) := \interpOp{\miset}[q](\y) := \sum_{\vMi \in \miset} c_{\vMi} \interpOp{\vMi}[q](\y),
\end{equation}
where the combination technique coefficients $c_{\vMi}$ and \ah{tensorisation of the interpolation operators \eqref{eq:oneDinterp}} are defined by
\begin{equation}
    c_{\vMi} := \sum_{\substack{\vMiAlt \in\{0,1\}^{\nY}\\ \vMi+\vMiAlt \in \miset}}(-1)^{\Vert\vMiAlt\Vert_1}, \quad \interpOp{\vMi}[q](\y):=\bigotimes_{\iY=1}^{\nY} \interpOp{\mi_{\iY}}[q](\y).
\end{equation}
Equation \eqref{eq:combinationtechnique} provides a global polynomial approximation over the parameter domain $\Gamma$, \mg{within the polynomial approximation space $\polySpace{\miset}(\Gamma;\R)$ determined by the multi-index set $\miset$.}
Denote the set of $\nColloc{\miset}$ distinct collocation points
\begin{equation} \label{eq:collocset}
    \mg{\collocSet{\miset}:=\bigcup_{\substack{\vMi\in\miset \\ c_{\vMi}\neq 0}} \collocSet{\mi_1} \times \collocSet{\mi_2} \times \cdots \times \collocSet{\mi_\nY}.}
\end{equation}
Note that if the sequence of points $\{\collocSet{\level}\}_{\level\in\N}$ is nested, i.e.\ $\collocSet{1}\subset \collocSet{2} \subset \collocSet{3} \subset \cdots$, then many points in \eqref{eq:collocset} are coincident.
\bk{The sparse grid construction \eqref{eq:collocset} is demonstrated with a $\nY=2$ dimensional example in \Cref{fig:sparsegridexample}.
    \begin{figure}
        \begin{subfigure}[t]{0.33\textwidth}
            \centering
            \inputfrom{figures/grids/}{tg1.tex}
            \caption{Tensor grid $\vMiParam=[1,2]$, $c_{\vMi}=-1$}
            \label{fig:tg1}
        \end{subfigure}
        \begin{subfigure}[t]{0.33\textwidth}
            \centering
            \inputfrom{figures/grids/}{tg2.tex}
            \caption{Tensor grid $\vMiParam=[1,3]$, $c_{\vMi}=1$}
            \label{fig:tg2}
        \end{subfigure}
        \begin{subfigure}[t]{0.33\textwidth}
            \centering
            \inputfrom{figures/grids/}{tg3.tex}
            \caption{Tensor grid $\vMiParam=[2,1]$, $c_{\vMi}=-1$}
            \label{fig:tg3}
        \end{subfigure}
        \begin{subfigure}[t]{0.33\textwidth}
            \centering
            \inputfrom{figures/grids/}{tg4.tex}
            \caption{Tensor grid $\vMiParam=[2,2]$, $c_{\vMi}=1$}
            \label{fig:tg4}
        \end{subfigure}
        \begin{subfigure}[t]{0.33\textwidth}
            \centering
            \inputfrom{figures/grids/}{tg5.tex}
            \caption{Tensor grid $\vMiParam=[3,1]$, $c_{\vMi}=1$}
            \label{fig:tg5}
        \end{subfigure}
        \begin{subfigure}[t]{0.33\textwidth}
            \centering
            \inputfrom{figures/grids/}{sparsegrid.tex}
            \caption{Sparse grid $\collocSet{\miset}$}
            \label{fig:sparsegrid}
        \end{subfigure}
        \caption{For a multi-index set $\miset=\{\vMiParam \in \N^2 \st \Vert \vMiParam \Vert_1 \leq 2 + 2\}$ using Clenshaw--Curtis points with the doubling rule from \eqref{eq:rules}, the tensor grids with non-zero combination technique coefficient $c_{\vMi}$ and the complete sparse grid $\collocSet{\miset}$ defined via \eqref{eq:collocset} are shown in \Cref{fig:tg1,fig:tg2,fig:tg3,fig:tg4,fig:tg5} and \Cref{fig:sparsegrid} respectively.}
        \label{fig:sparsegridexample}
    \end{figure}}

\paragraph{\ah{Admissibility of multi-index sets}}
A valid approximation in \eqref{eq:combinationtechnique} is obtained only if the multi-index set $\miset$ is admissible (also known as downwards closed, or a lower set).
Admissibility requires each multi-index $\vMi \in \miset$ to satisfy
\begin{equation}\label{eq:admissibile}
    \vMi - \vUnit_{\iY} \in \miset \forall \iY=1,...,\nY \st \mi_{\iY}>1
\end{equation}
where $\vUnit_{\iY}=[0,0,...,1,...,0]$ is the unit multi-index with \mg{value} $1$ only in the $i^\text{th}$ entry.
The set of multi-indices that can be added to $\miset$ whilst maintaining admissibility is called the reduced margin
\begin{equation}\label{eq:reducedmargin}
    \reducedmargin{\miset} := \Big\{\vMiParam\in\N^{\nY} \setminus \miset \st \vMiParam - \vUnit_{\iY} \in \miset \forall \iY=1,...,\nY \st \mi_{\iY}>1 \Big\}.
\end{equation}
For an admissible multi-index set $\miset\subset\N^{\nY}$ and nested set of points $\{\collocSet{\level}\}_{\level\in\N}$, the approximation \eqref{eq:combinationtechnique} is an interpolant.
\ah{An example multi-index set $\miset$ and the corresponding reduced margin $\reducedmargin{\miset}$ are shown in \Cref{fig:reducedmargin}.}
\begin{figure}\centering
    \inputfrom{figures/}{reducedmargin.tex}
    \caption{Example reduced margin $\reducedmargin{\miset}$ for a multi-index set with $\nY=2$.}
    \label{fig:reducedmargin}
\end{figure}

\paragraph{\ah{Polynomial approximation spaces}}
\mg{Recall that} $\polySpace{\miset}(\Gamma;\R)$ \mg{denotes} the polynomial space to which the approximation $q^{\miset}$ in \eqref{eq:combinationtechnique} belongs.
\ah{This space can be obtained as the span of a suitable set of multi-variate polynomials, for example the multi-variate monomials, or multi-variate spectral-type polynomials, e.g.\ Legendre or Chebysev polynomials.
We require a spectral-type polynomial basis, and again we denote by  $\{\oneChebBasis_{\miPoly}\}_{\miPoly\in\N}$ a sequence of univariate polynomials of degree $\miPoly$, and let $\chebBasis_{\vMiPoly}(\y)=\prod_{\iY=1}^{\nY} \bkrevise{\oneChebBasis}_{\miPoly_{\iY}}(y_\iY)$.}
For the polynomial approximation space $\polySpace{\miset}(\Gamma;\R)$ \lt{to which the sparse grid approximation \eqref{eq:combinationtechnique} belongs}, defined implicitly for multi-index set $\miset$, a corresponding multi-index set of polynomial degrees $\misetPoly(\miset)$ can be defined such that
\begin{equation}\label{eq:polsetPI}
    \text{ \bkrevisesout{$\{\chebBasis_{\vMiPoly}\}_{\vMiPoly\in\misetPoly(\miset)}$ where}{} } \polySpace{\miset}(\Gamma;\R) = \spanOp_{\vMiPoly\in\misetPoly(\miset)} \{\chebBasis_{\vMiPoly}\} ,
\end{equation}
i.e.\ $\misetPoly(\miset)$ is the appropriate multi-index set representing the multi-variate polynomial degrees.
For an approximation $\approxP{q}(\y)$ computed via the combination technique, it is straightforward to perform a change of basis to the multi-variate polynomial basis $\{\chebBasis_{\vMiPoly}\}_{\vMiPoly\in\misetPoly(\miset)}$ via the solution of a linear system \cite{Constantine2012,Formaggia2012}.
This gives a \emph{truncated} spectral polynomial expansion
\begin{equation} \label{eq:truncated}
    \approxP{q}(\y) \equiv \sum_{\vMiPoly\in\misetPoly(\miset)} \approxP{\FT{q}}_{\vMiPoly} \chebBasis_{\vMiPoly}(\y).
\end{equation}
The equivalence of the sparse grid formulation and the truncated spectral expansion \mg{is} at the core of \mg{the} \plateaumisc{} algorithm \mg{introduced in \Cref{sec:pmisc}}.
\mg{Indeed, once the MISC approximation is constructed, the change of basis to spectral polynomials in \eqref{eq:truncated} allows exploitation of the spectral polynomial coefficient decay properties in the constructed approximation.
    This information is thus employed to detect solver errors via the corruption of spectral polynomial coefficients (\Cref{sec:noisesc}) and consequently devise a strategy to construct an approximant robust to noisy evaluations (\Cref{sec:pmisc}).
}

\subsection{Effect of solver errors in sparse grid stochastic collocation}\label{sec:noisesc}
It is generally not possible to evaluate $q(\y)$ exactly: each evaluation used to build \eqref{eq:combinationtechnique} is an approximation $q^{\vMiModel}(\vCollocPt) \approx q(\vCollocPt)$
via a PDE solve where $\vMiModel$ \lt{is a multi-index of hyperparameters that control} the approximation fidelity (e.g.\ mesh size, solver tolerances,...).
\bk{Therefore,} a computable single-fidelity \ah{surrogate of the quantity of interest} can be formed using the combination technique \eqref{eq:combinationtechnique} with multi-index set $\miset$,
\begin{equation} \label{eq:approxcomputable}
    q(\y) \approx q^{\vMiModel,\miset}(\y) := \sum_{\vMi \in \miset} c_{\vMi} \interpOp{\vMi}[q^{\vMiModel}](\y).
\end{equation}
\bk{constructed using pointwise evaluations $q^{\vMiModel}(\vCollocPt)\approx q(\vCollocPt) \forall \vCollocPt \in \collocSet{\miset}$.}
Again, as in \eqref{eq:truncated}, a truncated spectral expansion can be \mg{obtained} via a change of basis,
\begin{equation}\label{eq:fullexpansion}
    q^{\vMiModel,\miset}(\y) \equiv \sum_{\vMiPoly\in\misetPoly(\miset)} \FT{q}^{\vMiModel,\miset}_{\vMiPoly} \chebBasis_{\vMiPoly}(\y).
\end{equation}
Each computable coefficient $\FT{q}^{\vMiModel,\miset}_{\vMiPoly}$ for $\vMiPoly\in\misetPoly(\miset)$ is an approximation \mg{of} the corresponding coefficient $\FT{q}_{\vMiPoly}$ from the \bk{full spectral expansion} \eqref{eq:spectral} and can be expressed as
\begin{equation} \label{eq:spectral-split}
    \FT{q}^{\vMiModel,\miset}_{\vMiPoly} = \underbrace{\FT{q}_{\vMiPoly}}_{(I)} + \underbrace{\FT{q}^{\miset}_{\vMiPoly} - \FT{q}_{\vMiPoly}}_{(II)} + \underbrace{\FT{q}^{\vMiModel,\miset}_{\vMiPoly} - \FT{q}^{\miset}_{\vMiPoly}}_{(III)}.
\end{equation}
The term $(I)$ in \eqref{eq:spectral-split} is the \emph{true} spectral series coefficient\mgrevise{, encapsulating information on the smoothness of the quantity of interest $q$. It follows that the decay of $(I)$ represents the theoretical behaviour of the spectral coefficients according to the regularity of $q(\y)$.} 
\mgrevise{In particular,} for an analytic function $q(\y)$ in a spectral-type polynomial expansion, we expect the expansion coefficients $\{\FT{q}_{\vMiPoly}\}_{\N_0^{\nY}}$ to obey an a priori bound of the form
\begin{equation} \label{eq:coeffdecay}
    \lvert \FT{q}_{\vMiPoly} \rvert \leq C(q,\vMiPoly) \prod_{\iY=1}^{\nY} e^{-g_{\iY}\miPoly_{\iY}}
\end{equation}
where $g_{\iY} > 0$ for all $\iY=1,2,...,\nY$ and $C(q,\vMiPoly)$ grows slowly as the polynomial degree $\vMiPoly$ increases.
For further details on analyticity and \bk{spectral polynomial expansion coefficient bounds} see, e.g.\, \cite{Nobile2008b}\cite[Chapter 3, Theorem 3.2]{Adcock2022}.
When the term $(I)$ dominates \eqref{eq:spectral-split} the coefficient is well approximated.

The term $(II)$ represents \bkrevisesout{the truncation error in the approximation}{}
\mgrevise{an approximation error due to computing the spectral coefficient $\FT{q}^{\miset}_{\vMiPoly}$ using the polynomial interpolant \eqref{eq:combinationtechnique} (rewritten in the form \eqref{eq:truncated}), instead of the true coefficient $\FT{q}_{\vMiPoly}$ in the infinite expansion \eqref{eq:spectral}.}
This contains aliasing type effects \mgrevise{due to using the subspace $\polySpace{\miset}(\Gamma;\R)$ in \eqref{eq:polsetPI} implicitly defined by the interpolation}.
This term is at its largest for the highest polynomial degree terms in the \bkrevisesout{truncated}{}approximation.
The effect of $(II)$ can be mitigated for a given polynomial degree $\vMiPoly$ by increasing the size of the polynomial approximation space.
\bkrevisesout{Truncation}{These} effects are generally not a source of concern in the context of this work.

The term $(III)$ represents the solver error\mgrevise{, which is usually considered as an intrinsic source of noise in simulation-based data. The proposed \plateaumisc{} strategy aims to robustify the standard MISC algorithm to handle this uncertainty.}
\bkrevise{By linearity} of the change of basis to the spectral polynomials and linearity of the combination technique construction \eqref{eq:combinationtechnique} \bkrevise{it follows that} $(III)$ is in fact the spectral coefficient of an interpolant of the solver error,
\begin{equation}\label{eq:interpoferror}
    \sum_{\vMiPoly\in\misetPoly(\miset)} \left(\FT{q}^{\miModel,\miset}_{\vMiPoly}- \FT{q}^{\miset}_{\vMiPoly}\right)  \chebBasis_{\vMiPoly}(\y)
    \equiv
    \bkrevise{q^{\vMiModel,\miset}(\y) - \approxP{q}(\y)\equiv}
    \sum_{\vMi \in \miset} c_{\vMi} \interpOp{\vMi}[\ah{q^{\alpha} - q}](\y).
\end{equation}
If the solver error is constant over $\Gamma$, i.e. \bkrevise{$E^{\vMiModel}(\y) = (q^{\vMiModel} - q)(\y) = C^{\vMiModel}\in\R$,}
then the interpolant on the right-hand side of \eqref{eq:interpoferror} is exact and a constant. 
Therefore \bkrevisesout{the term $(III)$ is a constant, and}{} the only non-zero coefficient \bkrevise{in the expansion on the left-hand side of \eqref{eq:interpoferror}} is for the constant polynomial, \bkrevisesout{(i.e.\ polynomial degree $\vMiPoly=\mistyle{0}$)}{that is, the term $(III)$ is non-zero only for $\vMiPoly=0$}.
If instead the solver error is more complex \bkrevisesout{and appears highly oscillatory}{}, possibly even seemingly random over $\Gamma$, many coefficients $\FT{q}^{\miModel,\miset}_{\vMiPoly}- \FT{q}^{\miset}_{\vMiPoly}$ will be non-zero and potentially large\bk{; \mgrevise{for instance,} this can be due to unresolved or preasymptotic meshes, coarse tolerances in iterative solvers, timestepping etc.}
It is this term $(III)$ that can corrupt the \mgrevise{computation of the} surrogate ${q}^{\vMiModel,\miset}$ in \eqref{eq:approxcomputable}.
This can result in a \emph{spectral plateau}, i.e.\ a situation in which many coefficients corresponding to sufficiently high polynomial degree $\Vert \miPoly\Vert_1$ terms in the expansion \eqref{eq:truncated} are corrupted and of a similar erroneous magnitude.
\bkrevise{The following Example \ref{ex:gaussian} provides further \mgrevise{intuition for} the expected stagnation in the decay of the spectral coefficients \mgrevise{under standard assumptions on the quantity of interest}.} 

\bkrevise{
  \begin{example}[Estimating coefficient plateau for Gaussian noise]\label{ex:gaussian}
    Suppose that the solver error $E^{\vMiModel}(\y)$ 
    is a white Gaussian random field over $\Gamma$, i.e., $E^{\vMiModel}(\y) \sim \mathcal{N}(0,\sigma^2)$ and
    $E^{\vMiModel}(\y_i)$ is statistically independent of $E^{\vMiModel}(\y_j)$ for all $\y_i, \y_j \in \Gamma$.
    \mgrevise{Consider the `linear' rule in \eqref{eq:rules}, $\nColloc{\level}=\level$, and take the abscissae of the mid-point quadrature rule in $[0,1]$
      to be the set of one-dimensional collocation points $\collocSet{\level} = \{ \collocPt_{\iColloc}\}_{\iColloc=1}^{\level}$, i.e., 
      $\collocPt_{\iColloc} = \frac{1}{\level} (\iColloc-\frac{1}{2})$; finally, let $\miset = [w,w,w,\cdots]$ for some $w \in \N, w\geq 1$.}
    The sparse grids interpolant reduces to a tensor-grid interpolant over a Cartesian grid
    $\collocSet{\miset} = \collocSet{w} \times \collocSet{w} \times \cdots=\{\vCollocPt_i\}_{i=1}^M$ with $M= w^d$ collocation points,
    and the coefficients $c_{\vMi}$ in \eqref{eq:combinationtechnique}, in \eqref{eq:approxcomputable} and in the right-hand side of \eqref{eq:interpoferror}
    are all zero other than $c_{[w,w,\cdots]}=1$.
    Furthermore, let $\{\chebBasis_{\vMiPoly_i}\}_{i=1}^M$ be the Legendre orthonormal polynomials with total polynomial degree up to $w-1$, i.e.
    $\misetPoly(\miset) = \{\vMiPoly : \|\vMiPoly\|_{\infty} \leq w-1\} $ in the left-hand side of \eqref{eq:interpoferror}.
    We introduce
    a suitable ordering of $\misetPoly(\miset)$ such that we can index its elements as $\vMiPoly_i$ and hence 
    the set of polynomials is $\{\chebBasis_{\vMiPoly_i}\}_{i=1}^M$ with a scalar index $i$.
    We can collect the pointwise errors $\{E^{\vMiModel}(\vCollocPt_i)\}_{i=1}^M$ in a (random) vector
    \begin{equation}
      \vec{e} :=
      \begin{bmatrix}
        E^{\vMiModel}(\vCollocPt_1),\, E^{\vMiModel}(\vCollocPt_2), \, \ldots,\, E^{\vMiModel}(\vCollocPt_M)
      \end{bmatrix}^{\top}
      \sim \mathcal{N}(\vec{0}, \sigma^2 \Id),
    \end{equation}
    where $\Id$ is the identity matrix in $\R^{M \times M}$.
    Similarly, we can collect the coefficients of the spectral expansion at the left-hand side of \eqref{eq:interpoferror}
    in another vector
    \[
      \FT{\vec{e}} :=
      \begin{bmatrix}
        \FT{q}^{\miModel,\miset}_{\vMiPoly_1} - \FT{q}^{\miset}_{\vMiPoly_1}, \, 
        \FT{q}^{\miModel,\miset}_{\vMiPoly_2} - \FT{q}^{\miset}_{\vMiPoly_2}, \,
        \ldots, \ldots
        \FT{q}^{\miModel,\miset}_{M} - \FT{q}^{\miset}_{M}
      \end{bmatrix}.
    \]
    The left-hand equality of \eqref{eq:interpoferror} can be written in compact form  as $\bm{Q} \FT{\vec{e}} = \vec{e}$, where
    \begin{equation}
      \bm{Q}:=\begin{bmatrix}
        \chebBasis_{\vMiPoly_1}(\vCollocPt_1) & \chebBasis_{\vMiPoly_2}(\vCollocPt_1) & \cdots & \chebBasis_{\vMiPoly_{\nColloc{\miset}}}(\vCollocPt_1)  \\
        \chebBasis_{\vMiPoly_1}(\vCollocPt_2) & \chebBasis_{\vMiPoly_2}(\vCollocPt_2) & \cdots & \chebBasis_{\vMiPoly_{\nColloc{\miset}}}(\vCollocPt_2)  \\
        \vdots                             & \vdots                              & \ddots & \vdots                                            \\
        \chebBasis_{\vMiPoly_1}(\vCollocPt_{M}) & \chebBasis_{\vMiPoly_2}(\vCollocPt_{M}) & \cdots & \chebBasis_{\vMiPoly_{\nColloc{\miset}}}(\vCollocPt_{M}) \\
      \end{bmatrix}.
    \end{equation}
    The spectral coefficients $\FT{\vec{e}} = \bm{Q}^{-1} \vec{e}$ are then jointly Gaussian, $\FT{\vec{e}} \sim \mathcal{N}(\vec{0},\Sigma)$
    with covariance matrix $\Sigma:= \sigma^2 \bm{Q}^{-1}(\bm{Q}^{-1})^{\top} = \sigma^2(\bm{Q}^{\top}\bm{Q})^{-1}$, where 
    each entry in the matrix $\bm{Q}^{\top}\bm{Q}$ can be written as
    \[
      \left[\bm{Q}^{\top}\bm{Q}\right]_{i,j} = \sum_{k=1}^{M} \chebBasis_{\vMiPoly_i}(\vCollocPt_k)\chebBasis_{\vMiPoly_j}(\vCollocPt_k).
    \]
    Next, we reinterpret $\frac{1}{M} \left[\bm{Q}^{\top}\bm{Q}\right]_{i,j}$ as a quadrature rule over $\Gamma$ for approximating
    $\left(\chebBasis_{\vMiPoly_i},\chebBasis_{\vMiPoly_j}\right)_{L^2(\Gamma)}$ and exploit the orthonormality of $\{\chebBasis_{\vMiPoly_i}\}_{i=1}^M$
    to obtain 
    \begin{equation}\label{eq:approxorthogonal}
      \frac{1}{M}  \left[\bm{Q}^{\top}\bm{Q}\right]_{i,j}
      = \frac{1}{M} \sum_{k=1}^{m_\miset} \chebBasis_{\vMiPoly_i}(\vCollocPt_k)\chebBasis_{\vMiPoly_j}(\vCollocPt_k)
      \approx \left(\chebBasis_{\vMiPoly_i},\chebBasis_{\vMiPoly_j}\right)_{L^2(\Gamma)}
      = \delta_{ij},
    \end{equation}
    where $\delta_{ij}$ is the classical Kronecker delta.
    This entails that $\bm{Q}^{\top}\bm{Q}$ is approximately diagonal, $\bm{Q}^{\top}\bm{Q} \approx M \Id$, and thus
    $\Sigma = \sigma^2(\bm{Q}^{\top}\bm{Q}  )^{-1} \approx \frac{1}{M} \sigma^2 \Id$.
    Therefore, the errors on the spectral coefficients are approximately statistically independent with variance $\FT{\sigma}^2 = \frac{1}{M} \sigma^2$.
    This implies that their absolute value approximately follows a folded normal distribution, and hence
    \begin{equation}
      \Exp{\left\lvert \FT{q}^{\miModel,\miset}_{\vMiPoly}- \FT{q}^{\miset}_{\vMiPoly} \right\rvert} = \FT{\sigma}\sqrt{2/\pi}
    \end{equation}
    for all $\vMiPoly\in\misetPoly(\miset)$.
    This results in a spectral plateau in the decay of $\FT{q}^{\miModel,\miset}_{\vMiPoly}$.
  \end{example}
}

\mgrevise{We will later observe a spectral plateau to be present in more general scenarios than that of Example \ref{ex:gaussian},
  cf. the numerical experiments of \Cref{sec:pmisc} and \Cref{sec:numerics}.
In \Cref{sec:pmisc} we present a method to detect the resulting spectral plateau effect, and embed this plateau detection into the adaptive MISC algorithm
to formulate an automatic strategy named \plateaumisc{}, that is robust to intrinsic solver noise.}

\section{PlateauMISC: \mg{a multi-fidelity algorithm robust to solver errors}} \label{sec:pmisc}
In this section we present the novel \plateaumisc{} algorithm.
This builds upon the classic MISC algorithm which is detailed in \Cref{sec:misc}.
A plateau detection algorithm is described in \Cref{sec:plateaudetection} and the fully \mg{automatic} \plateaumisc{} algorithm is presented in \Cref{sec:plateaumisc-algorithm}.
An illustrative example highlights the advantages of the proposed \plateaumisc{} algorithm in \Cref{example:2dgp}.

\subsection{Multi-Index Stochastic Collocation}\label{sec:misc}
Multi-index stochastic collocation (MISC) is a multi-fidelity extension of the sparse grid stochastic collocation combination technique \eqref{eq:combinationtechnique}.
The MISC algorithm has been investigated in \cite{Jakeman2019,HajiAli2016,HajiAli2016a,Zech2019,Beck2019} and applications have been investigated in \cite{Jakeman2022,Piazzola2022a,Chiappetta2024}.
Similar approaches (multi-level stochastic collocation) for parametric PDEs are considered in \cite{Teckentrup2015,Bespalov2022b}.

For a quantity of interest $q:\Gamma\to\R$, consider a hierarchy of computable \bkrevisesout{models}{fidelities} $q^{\vMiModel}(\y)\approx q(\y)$ indexed by multi-indices $\vMiModel \in \N^{\nModel}$.
Although a formal assumption is not made, it is expected that the \bkrevisesout{models}{fidelities} become more accurate as the \bkrevisesout{model}{} multi-index $\vMiModel$ increases in any dimension.
For example, one could consider a sequence of \bk{solvers} for an elliptic PDE problem given in domain $\Omega\subset \R^{\nModel}$ where mesh edge \mg{sizes} are proportional to $2^{-\miModel_{\iY}}$ for each spatial dimension $\iY=1,2,...,\nModel$.
Similarly, one could consider a time-dependent problem in which adaptive timestepping is applied with \mg{solvers using} sequentially smaller local error tolerances $\delta\propto 10^{-\miModel}$ as the \bkrevisesout{model}{fidelity} index $\miModel$ increases ($\nModel=1$ in this case).
\bkrevise{
We note that MISC does not directly apply when a hierarchical ordering of the fidelities cannot be established a priori (i.e., when increasing $\alpha$ does not necessarily lead to a reduction of the error in the evaluation of a quantity of interest). In such cases, it is necessary to employ methods that invest  some computational effort to identify correlations and effective hierarchies among the different fidelities; see, for example,
  \cite{gorodetsky.eal:ACV,narayan:bandit}.}

Consider an admissible multi-index set $\miset \subset \N^{\nModel+\nY}$ defined over both the \bkrevisesout{model}{fidelity} and parameter dimensions.
A multi-index $\vMiFull=[\vMiModel,\vMiParam]\in\miset$ splits with $\vMiModel\in\N^{\nModel}$ controlling the \bkrevisesout{model}{} fidelity $q^{\vMiModel}$ and $\vMiParam\in\N^{\nY}$ again defining an interpolation operator to give a term $\interpOp{\vMiParam}[q^{\vMiModel}]$.
The combination technique approximation \eqref{eq:combinationtechnique} is extended to give the MISC approximation,
\begin{equation}
    \label{eq:misc}
    q(\y)\approx\approxFinal{q}(\y) := \interpOp{\miset}[q](\y) := \sum_{\vMiFull \in \miset} c_{\vMiFull} \interpOp{\vMi}[q^{\vMiModel}](\y)  \text{ for all $\y\in\Gamma$ where } \vMiFull=[\vMiModel, \vMi].
\end{equation}
The resulting approximation is still a linear combination of global polynomial interpolants approximating the quantity of interest, with each interpolant acting on a different fidelity.
\ah{\begin{remark}[MISC is not interpolatory]
        Note that the resulting approximation \eqref{eq:misc} is not necessarily an interpolant: a collocation point $\vCollocPt$ may be requested by two interpolants acting on different fidelities $\vMiModel_1\neq\vMiModel_2$ with $q^{\vMiModel_1}(\vCollocPt)\neq q^{\vMiModel_2}(\vCollocPt)$, yet the approximation $\approxFinal{q}(\vCollocPt)$ can only take a single value.
    \end{remark}}

The greedy adaptive construction of the MISC approximation \eqref{eq:misc} is presented \mg{in} \Cref{alg:misc} \bk{\cite{Jakeman2019,HajiAli2016,Piazzola2022a}}, based upon the classical Gerstner--Griebel-type dimensional-adaptive loop \cite{Gerstner2003}. 
In particular, the construction follows \cite[Algorithm 1]{Piazzola2022a}.
\begin{algorithm}[tb]
    \caption{Multi-Index Stochastic Collocation}
    \label{alg:misc}
    \begin{algorithmic}[1]
        \Function{MISC}{$\{q^{\vMiModel}:\Gamma\subset\R^{\nY} \to \R\}_{\vMiModel \in \N^{\nModel}}$}
        \State Initial multi-index set $\miset = \{[\mistyle{1},\mistyle{1}]\}\subset \N^{\nModel+\nY}$.
        \While{stopping criteria not met}
        \State{Compute approximation $\approxFinal{q}:=\sum_{\vMiFull\in\miset} c_{\vMiFull} \interpOp{\vMi}[q^{\vMiModel}]$ where $\vMiFull=[\vMiModel, \vMiParam]$.} \Comment{Solve}
        \State{Compute reduced margin $\reducedmargin{\miset}$} via Equation \eqref{eq:reducedmargin}.
        \For{$\vMiFull \in \reducedmargin{\miset}$} \Comment{Estimate}
        \State Estimate multi-index contribution $E_{\vMiFull}$ via Equation \eqref{eq:error-ind}.
        \State Estimate multi-index cost $W_{\vMiFull}$ via Equation \eqref{eq:work-misc}.
        \State Compute multi-index profit $\profit_{\vMiFull}=E_{\vMiFull} / W_{\vMiFull}$.
        \EndFor

        \State Select $\vMiFull^* \gets \argmax_{\vMiFull\in\reducedmargin{\miset}} \profit_{\vMiFull}$ and update $\miset \gets \miset \cup \vMiFull^*$. \Comment{Mark and Refine}
        \EndWhile
        \State{\Return $\approxFinal{q}:=\sum_{\vMiFull\in\miset\cup\reducedmargin{\miset}} c_{\vMiFull} \interpOp{\vMi}[q^{\vMiModel}]$ where $\vMiFull=[\vMiModel, \vMiParam]$.}
        \EndFunction
    \end{algorithmic}
\end{algorithm}
\Cref{alg:misc} is applied using a \mg{set of solvers giving approximations} $\{q^{\vMiModel}\}_{\vMiModel\in\N^{\nModel}}$ where $\vMiModel$ defines the fidelity.
It iteratively grows the multi-index set $\miset$ and consequently grows the polynomial approximation space $\polySpace{\miset}(\Gamma;\R)$.
Loosely speaking, adding a multi-index $\vMiFull=[\vMiModel, \vMiParam]$ in which we increase $\vMiModel$ refines the fidelity \bkrevisesout{of the model}{} used by an interpolant defined via $\vMiParam$, whilst increasing $\vMiParam$ will increase the number of evaluations and parametric approximation space for a particular fidelity $\vMiModel$.

At each iteration, the neighbouring multi-indices in the reduced margin $\reducedmargin{\miset}$ \eqref{eq:reducedmargin} are considered for refining the current multi-index set $\miset$.
Error indicators $ E_{\vMiFull}$ are computed for each $\vMiFull\in\reducedmargin{\miset}$.
\mg{To define these, first it is assumed that for each parameter $y_\iY$ there is a corresponding weight function $\rho_{\iY}:\Gamma\to\R$ for $\iY=1,2,...,\nY$.
    For simplicity these are assumed to be $\rho_{\iY}\equiv 1$ for $\iY=1,2,...,\nY$, i.e.\ each parameter can be considered to be representative of the image of an independent, uniformly distributed random variable on $[0,1]$.
    Extension to other commonly used parameter domains and associated weight functions is straightforward \lt{and well established in the UQ literature \cite{ghanem:uqbook,Smith2013,Sullivan2015}}}.
The error indicators are defined as
\begin{equation} \label{eq:error-ind}
    \begin{aligned}
        E_{\vMiFull} & := \Big\lvert \Exp{\interpOp{\miset \cup \{\vMiFull\}}[{q}](\y)} - \Exp{\interpOp{\miset}[q](\y)} \Big\rvert                                                                                                                            \\
                     & := \Bigg\lvert \int_{\Gamma}{\interpOp{\miset \cup \{\vMiFull\}}[{q}]}(\y) \rho(\y)  \dd \y - \int_{\Gamma}{\interpOp{\miset}[q]}(\y) \rho(\y)\dd \y \Bigg\rvert \quad \text{ where } \rho(\y) = \prod_{\iY=1}^{\nY} \rho_{\iY}(y_\iY).
    \end{aligned}
\end{equation}
\mg{Equation \eqref{eq:error-ind}} characterises the change in the approximation due to the addition of a multi-index $\vMiFull$.
Alternative estimators can be used, for example considering the \lt{$L_{\rho}^2(\Gamma;\R)$ or $L_{\rho}^\infty(\Gamma;\R)$} norm of the difference $\interpOp{\miset \cup \{\vMiFull\}}[{q}] - \interpOp{\miset}[q]$.
Next, for each fidelity $\vMiModel$, a representative fidelity solve cost $\widehat{W}_{\vMiModel}$ for a single approximation \mg{is introduced}.
The cost for adding a multi-index $\vMiAlt=[\vMiModel, \vMiParam]$ to the multi-index set $\miset$ is the number of additional solves $\nColloc{\miset\cup\{\vMiFull\}}\rvert_{\vMiModel} - \nColloc{\miset}\rvert_{\vMiModel}$ (i.e.,\ the number of additional sparse grid points using fidelity $\miModel$\bk{, assuming nestedness of the one-dimensional collocation points $\{\collocSet{\level}\}_{\level\in\N}$}) multiplied by the representative fidelity solve cost $\widehat{W}_{\vMiModel}$, \mg{that is}
\begin{equation} \label{eq:work-misc}
    W_{\vMiFull} := \widehat{W}_{\vMiModel}(\nColloc{\miset\cup\{\vMiFull\}}\rvert_{\vMiModel} - \nColloc{\miset}\rvert_{\vMiModel}).
\end{equation}
The error indicator and \mg{cost} are then combined to give a profit $\profit_{\vMiFull}$ for each $\vMiFull\in\reducedmargin{\miset}$, \mg{namely,}
\begin{equation}
    \profit_{\vMiFull}:= E_{\vMiFull}/ W_{\vMiFull}.
\end{equation}
The greedy method in \Cref{alg:misc} selects the largest profit and refines the approximation using the corresponding multi-index $\vMiFull^*$.
The process repeats until a stopping criterion, for example a maximum total cost, maximum number of solves or a \bk{threshold} profit, is reached.
The final returned approximation uses the multi-index set $\miset \cup \reducedmargin{\miset}$.
This has no additional cost because the evaluations required to include $\reducedmargin{\miset}$ have already been performed when computing the error indicators.


\subsection{Effect of solver error in MISC}
\ah{
    The solver error discussed in \Cref{sec:noisesc} can lead the greedy \mg{approach} in \Cref{alg:misc} to mistakenly prioritise ``noisy'' fidelities during the refinement process.
    Each refinement step effectively adds higher-degree polynomial terms to the approximation space $\polySpace{\miset}(\Gamma;\R)$.
    In the interpolation basis, these higher-degree terms may begin to fit the solver error rather than the true solution;
    equivalently, in the spectral polynomial basis, this refinement may introduce spectral terms with spurious, large coefficients, particularly in the spectral plateau, corresponding to the term $(III)$ in the splitting \eqref{eq:spectral-split}.
    This leads to an inflated error indicator \eqref{eq:error-ind}, causing the greedy \Cref{alg:misc} to erroneously ``chase the noise''\mg{: it will tend to overfit the solver error, being unable to identify the \emph{true} solution behaviour}.
    As a result, the adaptive procedure may overuse inaccurate, noisy (low fidelity) evaluations with falsely large error indicators, while overlooking more accurate fidelities.
    This behaviour can corrupt the approximation, resulting in a noisy and highly oscillatory surrogate $q^{MISC}(\y)$.
    \mg{An example of this phenomenon will be discussed later in \Cref{example:2dgp}}.
}

\subsection{Plateau detection} \label{sec:plateaudetection}

\bk{To detect solver noise in the MISC approximation, we introduce a plateau detection algorithm}.
It exploits the truncated spectral polynomial expansion \eqref{eq:fullexpansion} and the expected spectral coefficient decay \eqref{eq:coeffdecay} to detect the effect of solver error via the spectral coefficient plateau resulting from term $(III)$ in \eqref{eq:spectral-split}. 
The plateau detection borrows ideas from the approximation package \chebfun{} and its \texttt{standardChop} algorithm \cite[Table 3.3]{Boyd2014}\cite{Aurentz2017}.
\bk{Similar ideas have also been used in the context of shock capturing, where elementwise spectral-type polynomial expansions are used to detect non-smooth solutions \cite{Persson2006,Huerta2011}}.

We begin our discussion with the single-fidelity case \ah{before extending to} the multi-fidelity \ah{setting}.
Consider the expansion \eqref{eq:fullexpansion} of the single fidelity approximation $q^{\vMiModel,\miset}(\y)$ for a \ah{given} fidelity \ah{level} $\vMiModel$.
For an approximation using the multi-index set $\miset$, there is a corresponding set of polynomial degrees $\misetPoly(\miset)$ and corresponding computable spectral polynomial coefficients $\{\FT{q}^{\vMiModel,\miset}_{\vMiPoly}\}_{\vMiPoly \in \misetPoly(\miset)}$ \bk{defining the truncated spectral polynomial expansion}.
\ah{Let} $k_{e}:=\max_{\vMiPoly\in\misetPoly({\miset})} \Vert \vMiPoly\Vert_1$ \ah{be} the maximum polynomial total degree in the \bk{truncated} spectral polynomial expansion.
\ah{We define a} \emph{spectral envelope} $\envelope:\envI\subset\N_0 \to \R$ \ah{with} $\envI:=\{0,1,...,k_{e}\}$ \ah{to provide an upper} bound \ah{for} the \ah{magnitude of the} spectral coefficients, \ah{specifically}
\begin{equation} \label{eq:env}
    \envelope(i):=\max_{\vMiPoly \in \misetPoly({\miset}) \st \Vert \vMiPoly \Vert_1 \geq i}\big\lvert \FT{q}^{\vMiModel,\miset}_{\vMiPoly}\big\rvert.
\end{equation}
The envelope $\envelope$ is expected to show two regimes: an exponential decrease as predicted by the \bk{a priori bound }\eqref{eq:coeffdecay}, followed by a spectral plateau due to solver error as predicted by $(III)$ in \eqref{eq:spectral-split}.
This \ah{behaviour} will be illustrated in the forthcoming numerical examples.
\ah{To identify this plateau, we fit} a piecewise linear model \ah{to} $\log_{10}(\envelope)$ and \ah{detect flattening in the tail.}
Note that the \ah{a priori decay rate} \eqref{eq:coeffdecay} \ah{may not accurately describe the low polynomial degree terms.}
\ah{Therefore, we ignore the initial $\nBurnIn$ entries of the envelope (\emph{burn in})}.
Similarly, the highest polynomial degree spectral coefficients may be poorly estimated \ah{due to} aliasing  effects (i.e., regime (II) in splitting \eqref{eq:spectral-split}) \ah{so we discard the} final $\nBurnOut$ terms (\emph{burn out}) \ah{in the plateau detection process}.

In practice, a \emph{change point detection} algorithm with a piecewise linear model is used to partition $\{\left(i,\log_{10}(\envelope(i))\right)\}_{i=\nBurnIn}^{k_e - \nBurnOut}$ into two with the change point $\changept$ \cite{Killick2012}.
These fitted functions will later be referred to as \emph{piecewise log-linear} functions.
The models are defined by minimising the sum of the square errors,
\begin{equation}\label{eq:findchangepts}
    \begin{aligned}
        \left\{\changept, k_0, c_0, k_1, c_1\right\} = \textrm{argmin}_{\substack{j\in\{\nBurnIn,...,k_e - \nBurnOut\} \\
        m_0, c_0, m_1, c_1 \in \R}} \sum_{i=\nBurnIn}^{j-1} & \left(\log_{10}(e(i)) - m_0 i - c_0\right)^2             \\&+ \sum_{i=j}^{k_e-\nBurnOut} \left(\log_{10}(e(i)) - m_1 i - c_1\right)^2.
    \end{aligned}
\end{equation}
The change point $\changept$ is returned along with two linear models
\begin{equation}\label{eq:linearmodels}
    \begin{cases}
        m_0 i + c_0 , & i \in \{\nBurnIn,...,\changept-1\},        \\
        m_1 i + c_1 , & i \in \{\changept,...,k_{e} - \nBurnOut\}. \\
    \end{cases}
\end{equation}
The second model is associated with the spectral plateau.
As already hinted, our algorithm identifies a plateau if the following conditions are met:
\begin{enumerate}
    \item the second linear model has a small gradient $\lvert m_1 \rvert \leq m_*$,
    \item the second model is suitably long $k_{e} - \nBurnOut - \changept > \nMin$.
\end{enumerate}
The spectral plateau level is computed as $\error{plateau} := 10^{m_1\changept +c_1}$.
Spectral plateau detection is summarised \mg{in} \Cref{alg:plateau}.
In practice, we use the optimised \matlab{} \texttt{findchangepts} solver to identify the models in \eqref{eq:findchangepts} \mg{and \eqref{eq:linearmodels}}.
Although ideas are drawn from the \texttt{standardChop} algorithm in \cite{Aurentz2017}, the  \texttt{standardChop} algorithm is not appropriate in our context because it is designed for identification of spectral plateau effects at or near numerical precision.

\bkrevise{Four hyperparameters have been introduced in the plateau detection algorithm
    detailed above ($\nBurnIn$, $\nBurnOut$, $m_*$, $\nMin$). The following considerations can be useful when tuning them:
    \begin{itemize}
        \item Burn in $\nBurnIn$: This should be tuned based upon the expected decay of the spectral coefficients of the target response surface, $\FT{q}_{\vMiPoly}$, in the ideal case in which the solver noise is absent (see Eq. \ref{eq:spectral-split}). If the surface is expected to be extremely smooth it can be set equal to zero. If conversely the surface is expected to have e.g. local features, steep gradients, etc. there will be a pre-asymptotic behaviour that must be captured by the spectral polynomial expansion before the exponential decay in coefficient magnitude is observed. From our experiments, we suggest a burn in of $\nBurnIn = 2$.
        \item Burn out $\nBurnOut$: This is tuned to remove any aliasing effects in the tail of the spectral envelope. \mgrevise{An example of such effects will be presented in the numerical results of \Cref{sec:parabolic}}. From our numerical experiments, we suggest a burn out of $\nBurnOut=2$.
        \item Gradient threshold $m_*$: A gradient threshold of zero is almost impossible to attain in practice, whilst a gradient of $O(1)$ would still indicate spectral coefficient decay.
        We suggest a gradient threshold of $m_*=0.1$, which could be decreased if the response surface is expected \mgrevise{to experience} a fast coefficient decay.
        \item Minimum plateau length $\nMin$: The minimum plateau length must be chosen to balance the benefit of ensuring all useful information is extracted from cheaper, \mgrevise{low fidelity} simulations, whilst avoiding \mgrevise{unnecessary} computational \mgrevise{burden} oversampling the spectral plateau. We suggest a minimum plateau length $\nMin=3$. This could be increased if higher fidelity solvers are significantly more expensive, making it crucial to extract as much information as possible from the cheaper solvers.
    \end{itemize}}
    \MGnew{It is worth noting that the hyperparameter values are problem-dependent. Nonetheless, the PlateauMISC strategy exhibited behaviour that was robust to their selection. The numerical experiments in \Cref{sec:pmisc} and \Cref{sec:numerics} followed these principles to determine suitable values.}
 
\begin{algorithm}[tb]
    \caption{Spectral plateau detection}
    \label{alg:plateau}
    \begin{algorithmic}[1]
        \Function{PlateauDetection}{$\{i, \envelope(i)\}_{\bkrevise{i=0}}^{\bkrevise{k_e}}, \nBurnIn, \nBurnOut, \nMin$}
        \State{Compute $[\changept, m_0, c_0, m_1, c_1]$ via \eqref{eq:findchangepts}}


        \If{$\left|m_1\right| < m_*$ and $k_{e} - \nBurnOut - \changept > \nMin$}
        \State $\error{plateau} \gets 10^{m_1\changept +c_1}$
        \State ${isPlateau} \gets \texttt{true}$
        \Else
        \State $\error{plateau} \gets 0$
        \State ${isPlateau} \gets \texttt{false}$
        \EndIf

        \State \textbf{return} $\changept, \error{plateau}, {isPlateau}$
        \EndFunction
    \end{algorithmic}
\end{algorithm}

\subsection{\plateaumisc{} variant of MISC}\label{sec:plateaumisc-algorithm}
The plateau detection \ah{procedure described in} \Cref{alg:plateau} \mg{is} now integrated into the adaptive MISC approximation \ah{framework} of \Cref{alg:misc}.
\ah{To achieve this,} \Cref{alg:misc} is modified by \ah{incorporating} a \texttt{detect} step in the \ah{standard} \texttt{solve-estimate-mark-refine} loop, \ah{enabling the identification of} the spectral plateaus.
Since MISC is a multi-fidelity approximation method, plateau detection \ah{is applied independently at each fidelity level}.

Consider an approximation $\approxFinal{q}$ of the form \eqref{eq:misc}.
Denote the set of active fidelities in the approximation by $\modelSet$, that is,
\begin{equation}\label{eq:activemodels}
    \modelSet = \{\vMiModel \text{ such that there exists a $\vMi$ for which } [\vMiModel,\vMi]\in\miset\}.
\end{equation}
For each active fidelity $\vMiModel\in\modelSet$, a restriction to the fidelity is formed:
\begin{equation}\label{eq:miscrestricted}
    \approxFinal{q}\rvert_{\vMiModel}:=\interpOp{\miset{\rvert_{\vMiModel}}}[q] \text{ where } \miset\rvert_{\vMiModel}:={} \left\{ \widehat{\vMiParam} \st [\widehat{\vMiModel}, \widehat{\vMiParam}] \in \miset \st \widehat{\vMiModel} = \vMiModel \right\}.
\end{equation}
\bk{Each multi-index set $\miset{\rvert_{\vMiModel}}$ is downward-closed by construction, therefore \eqref{eq:miscrestricted} defines valid single-fidelity SGSC approximations.}
For each restricted approximation $\approxFinal{q}\rvert_{\vMiModel}$, an envelope $\envelope\rvert_{\vMiModel}$ can be defined via \eqref{eq:env}, and change points $\changept\rvert_{\vMiModel}$ and plateau levels $\error{plateau}\rvert_{\vMiModel}$ can be computed using \Cref{alg:plateau}.
Fidelities that have been identified as having a spectral plateau are collected as the set of \emph{saturated fidelities} $\saturatedModels\subset\N^{\nModel}$.
High polynomial degree terms for saturated fidelities should not be marked for further parametric refinement, \mg{limiting their employment in the construction of the surrogate model.}





To stop the greedy algorithm \mg{from} adding high polynomial degree terms for the saturated fidelities $\saturatedModels$, \mg{the \texttt{mark-refine} step in \Cref{alg:misc}} is preceded by setting a subset of profits to zero.
This subset corresponds to multi-indices $\vMiFull=[\vMiModel,\vMiParam] \in \reducedmargin{\miset}$ that \bk{satisfy the following criteria (see also \Cref{fig:filterprofits}):}
\begin{itemize}
    \item $\vMiFull$ corresponds to a saturated fidelity $\vMiModel\in\saturatedModels$;
    \item $\vMiFull$ only adds polynomial terms $ \vMiPoly \in \misetPoly(\miset\rvert_{\vMiModel} \cup\vMiParam) \setminus \misetPoly(\miset\rvert_{\vMiModel})$ of polynomial total degree greater than or equal to the fidelity's change point $\changept\rvert_{\vMiModel}$, that is if $\Vert \vMiPoly \Vert_1 \geq \changept\rvert_{\vMiModel}$ for all $ \vMiPoly \in \misetPoly(\miset\rvert_{\vMiModel} \cup\vMiParam) \setminus \misetPoly(\miset\rvert_{\vMiModel})$.
\end{itemize}
The second point allows low polynomial total degree terms to still be added via refinement, even for saturated fidelities.
These conditions can be expressed as
\begin{equation}\label{eq:filter-profits}
    \profit_{\vMiFull} \gets 0 \text{ if }\vMiModel \in \saturatedModels \text{ and } \Vert \vMiPoly \Vert_1 \geq \changept\rvert_{\vMiModel} \forall \vMiPoly \in \misetPoly(\miset\rvert_{\alpha} \cup\vMiParam) \setminus \misetPoly(\miset\rvert_{\alpha})
\end{equation}
where $\vMiFull=[\vMiModel,\vMi]$ and $\changept\rvert_{\vMiModel}$ is the change point for the associated fidelity $\vMiModel$ computed by \Cref{alg:plateau}.

\bk{\begin{figure}
        \begin{subfigure}[t]{0.49\textwidth}
            \inputfrom{figures/miset/}{miset_filtering.tex}
            \caption{\bk{Multi-index set $\miset$, reduced margin $\reducedmargin{\miset}$ and subset of reduced margin with profits $\profit_{\vMiParam}$ set to zero.}}
            \label{fig:miset-filtering}
        \end{subfigure}
        \begin{subfigure}[t]{0.49\textwidth}
            \inputfrom{figures/miset/}{polyspace.tex}
            \caption{\bk{Polynomial space $\misetPoly(\miset)$ corresponding to the multi-index set $\miset$, and the additional polynomials $\misetPoly(\miset\cup\vMiParam)\setminus\misetPoly(\miset)$ added by the multi-indices $\vMiParam \in \reducedmargin{\miset}$.}}
            \label{fig:polyspace}
        \end{subfigure}
        \caption{\bk{Example demonstrating the filtering of profits specified in \eqref{eq:filter-profits}.
                An example multi-index set $\miset$ is shown in \Cref{fig:miset-filtering}.
                The polynomial space $\misetPoly(\miset)$ is shown in \Cref{fig:polyspace} when using $\miset$ with the `linear' level-to-knots rule from \eqref{eq:rules}.
                If the corresponding fidelity is saturated with change point $\changept=3$, the profits for multi-indices $[4,2]$ and $[5,1]$ are set to zero because they only add polynomial terms of total degree $\Vert \vMiPoly\Vert_1 \geq 3$.
                We still permit the multi-index $[1,3]$ because it adds at least one polynomial term of total degree $\Vert \vMiPoly\Vert_1 < 3$.}}
        \label{fig:filterprofits}
    \end{figure}}

\ah{
    Another technical challenge arises from the treatment of the {saturated fidelities} $\saturatedModels$.
    Blocking further refinement of these saturated fidelities may unintentionally restrict access to high-degree polynomial approximations at \emph{unsaturated} fidelity levels.}
This \ah{issue stems from} the admissibility \ah{constraint} on the multi-index set, which is enforced by only exploring the reduced margin $\reducedmargin{\miset}$ \ah{as defined in equation} \eqref{eq:reducedmargin}.
For instance, \ah{consider} \Cref{fig:modifiedreducedmargin}\lt{:} the multi-index $\vMiFull=[3,3]$ may not be reached if the fidelities $\miModel=1,2$ are already saturated.
However, fidelity $\miModel=3$ \ah{being} less noisy, may \ah{still benefit from further parametric refinement}.
\begin{figure}\begin{center}
        \inputfrom{figures/}{modifiedreducedmargin.tex}
    \end{center}
    \caption{Modified reduced margin $\reducedmargin{\miset,\saturatedModels}$ defined in \eqref{eq:modifiedreducedmargin} for a multi-index set with $\nModel=1$ and $\nY=1$.
        If the saturated \bkrevisesout{models}{fidelities} are $\saturatedModels=\{1,2\}$, the multi-index $[3,3]$ is still permitted to be explored by the adaptive algorithm, whilst the multi-indices $[1,3],[2,3]$ may have their profits set to zero by \eqref{eq:filter-profits}.
        To maintain admissibility the multi-indices $B_{\miset,[3,3]}$, as defined in \eqref{eq:backfill}, must be used to `backfill` the multi-index set if $[3,3]$ is selected by the adaptive algorithm.}
    \label{fig:modifiedreducedmargin}
\end{figure}
\ah{To enable this}, we \ah{introduce} a \emph{modified reduced margin} \ah{that permits selective refinement at unsaturated fidelities}\bk{:
    first define the \bk{margin} set \ah{$\margin{\miset}$}},
\begin{equation}\label{eq:margin}
    \margin{\miset{}} := \Big\{\mistyle{\mu}\in\N^{\nModel+\nY} \st \mistyle{\mu} = \vMiFull + \vUnit_{\iY} \text{ for some } \{\vMiFull, \iY\} \in \miset \times \{1,2,...,\nModel+\nY\}\Big\}
\end{equation}
then the \emph{modified reduced margin} is
\begin{equation}\label{eq:modifiedreducedmargin}
    \reducedmargin{\miset,\saturatedModels} := \Big\{ \vMiFull \in \margin{\miset} \st [\vMiModel, \vMiParam]:= \vMiFull - \vUnit_{\iY} \in \miset \forall \iY=1,2,...,\nModel+\nY \text{ or } \vMiModel \in \saturatedModels\}.
\end{equation}
\ah{This construction ensures that each} backward neighbour $\vMiAlt - \vUnit_{\iY} $ either \ah{belongs to} the \ah{current} multi-index set $\miset$ or corresponds to a saturated fidelity in $\saturatedModels$ (i.e.\ $\vMiAlt - \vUnit_{\iY}\in \bigcup_{\vMiModel\in\saturatedModels} \{[\vMiModel, \vMiParam] \forall \vMiParam \in \N^{\nY}\}$).
\bk{Multi-indices corresponding to high polynomial degree terms in saturated fidelities are not explicitly removed from the reduced margin, however they \mg{are} ignored in the adaptive algorithm \lt{since} their profits are set to zero by \eqref{eq:filter-profits}.}

\ah{T}o maintain admissibility of the multi-index set when adding a multi-index from the modified reduced margin $\reducedmargin{\miset,\saturatedModels}$, we introduce a backfilling procedure.
For $\vMiFull\in \reducedmargin{\miset,\saturatedModels}$ the indices required to backfill the multi-index set $\miset$ to maintain admissibility are denoted by
\begin{equation}\label{eq:backfill}
    B_{\miset,\vMiFull} := \{\text{minimal set of } \widehat{\vMiFull}\in\N^{\nModel+\nY} \st \miset \cup B_{\miset,\vMiFull} \cup \{\vMiFull\} \text{ satisfies admissibility property \eqref{eq:admissibile}} \},
\end{equation}
see again \Cref{fig:modifiedreducedmargin} for a visualisation.
By construction all multi-indices $\vMiFull\in B_{\miset,\vMiFull}$ can only correspond to saturated fidelities.
The refinement of the multi-index set $\miset$ with $\vMiFull\in\reducedmargin{\miset,\saturatedModels}$ \ah{is then updated as} $\miset \gets \miset \cup B_{\miset,\vMiFull} \cup \{\vMiFull\}$.
Finally, the introduction of backfilling requires that we update the definitions of the error indicator $E_{\vMiFull}$ and cost $W_{\vMiFull}$ for a multi-index $\vMiFull\in\reducedmargin{\miset,\saturatedModels}$ accounting for the backfill set $B_{\miset,\vMiFull}$.
The error indicator \eqref{eq:error-ind} is \mg{thus} updated to
\begin{equation} \label{eq:error-ind-pmisc}
    E_{\vMiFull} := \Big\lvert \Exp{\interpOp{\miset \cup B_{\miset,\vMiFull} \cup  \{\vMiFull\}}[{q}]} - \Exp{\interpOp{\miset}[q]} \Big\rvert,
\end{equation}
and similarly the cost \eqref{eq:work-misc} is updated to
\begin{equation} \label{eq:work-pmisc}
    W_{\vMiFull} := \widehat{W}_{\vMiModel}(\nColloc{\miset\rvert_{\vMiModel} \cup\{\vMiParam\}} - \nColloc{\miset\rvert_{\vMiModel} }) + \sum_{\widehat{\vMiFull}=[\widehat{\vMiModel},\widehat{\vMiParam}]\in B_{\miset,\vMiFull}} \widehat{W}_{\widehat{\vMiModel}} (\nColloc{\miset\rvert_{\widehat{\vMiModel}} \cup\{\widehat{\vMiParam}\}} - \nColloc{\miset\rvert_{\widehat{\vMiModel} }})
\end{equation}
where the cost required to backfill the multi-index set must be included.

\begin{algorithm}[tb]
    \caption{Plateau Detection Multi-Index Stochastic Collocation}
    \label{alg:pmisc}
    \begin{algorithmic}[1]
        \Function{PlateauMISC}{$\{q^{\vMiModel}:\Gamma\subset\R^{\nY} \to \R\}_{\vMiModel \in \N^{\nModel}}$}
        \State Initial multi-index set $\miset = \{[\mistyle{1},\mistyle{1}]\}\subset \N^{\nModel+\nY}$
        \While{stopping criteria not met}
        \State{Compute approximation $\approxFinal{q}:=\sum_{\vMiFull\in\miset} c_{\vMiFull} \interpOp{\vMi}[q^{\vMiModel}]$ where $\vMiFull=[\vMiModel, \vMiParam]$.} \Comment{Solve}
        \State{Compute modified reduced margin $\reducedmargin{\miset,\saturatedModels}$ via Equation \eqref{eq:modifiedreducedmargin}.}
        \For{$\vMiFull \in \reducedmargin{\miset,\saturatedModels}$} \Comment{Estimate}
        \State Estimate multi-index contribution $E_{\vMiFull}$ via Equation \eqref{eq:error-ind-pmisc}.
        \State Estimate multi-index cost $W_{\vMiFull}$ via Equation \eqref{eq:work-pmisc}.
        \State Compute multi-index profits $\profit_{\vMiFull}=E_{\vMiFull} / W_{\vMiFull}$.
        \EndFor
        \State{}
        \For{each unsaturated active  \bkrevisesout{model}{fidelity} $\vMiModel \in \modelSet\setminus\saturatedModels$} \Comment{Detect}
        \State Construct $\approxFinal{q}\rvert_{\vMiModel}$ via \eqref{eq:miscrestricted}. 
        \State Construct spectral envelope $(i, \envelope\rvert_{\vMiModel}(i))$ via \eqref{eq:env}.
        \State $\changept\rvert_{\vMiModel}, \error{plateau}\rvert_{\vMiModel}, {isPlateau}\rvert_{\vMiModel} = \textsc{PlateauDetection}{(i, \envelope\rvert_{\vMiModel}(i), \nBurnIn, \nBurnOut, \nMin)}$
        \EndFor
        \State Update saturated  \bkrevisesout{models}{fidelities} $\saturatedModels \gets \saturatedModels\cup\{\vMiModel\in \modelSet\setminus\saturatedModels \st {isPlateau}\rvert_{\vMiModel} = \texttt{true}\}$
        \State{}
        \State \mg{Set} zero profits $\profit_{\vMiFull}$ for saturated $\vMiFull$ using criteria in \eqref{eq:filter-profits}.
        \State Select $\vMiFull^* \gets \argmax_{\vMiFull\in\reducedmargin{\miset,\saturatedModels}} \profit_{\vMiFull}$ and refine $\miset \gets \miset \cup B_{\miset,\vMiAlt^*} \cup\{\vMiFull^*\}$. \Comment{Mark and Refine}
        \EndWhile
        \State{\Return $\approxFinal{q}:=\sum_{\vMiFull\in\miset} c_{\vMiFull} \interpOp{\vMi}[q^{\vMiModel}]$ where $\vMiFull=[\vMiModel, \vMiParam]$.}
        \EndFunction
    \end{algorithmic}
\end{algorithm}
The complete \plateaumisc{} algorithm is presented in \Cref{alg:pmisc}.
Note the differences \mg{with respect} to \Cref{alg:misc}:
\begin{itemize}
    \item the \mg{plateau} detection step is introduced to identify saturated fidelities;
    \item a subset of profits may be set to zero via \eqref{eq:filter-profits};
    \item the returned approximation no longer uses the reduced margin.
\end{itemize}
This third point is to ensure that the high degree polynomial terms used to identify the spectral plateaus in the saturated fidelities are not included in the returned approximation, \bk{since they} will generally only contribute spurious oscillations to the \bk{full} approximation.

\subsection{Illustrative validation example}
\label{example:2dgp}
The \plateaumisc{} algorithm is validated using an analytic two-dimensional Gaussian peak (2DGP) test function from the Genz suite of test functions \cite{Genz1984}.
Consider the function $q:\Gamma:=[0,1]^2\to\R$
\begin{equation} \label{eq:2dgp}
    q(\y):=  \exp(-C_1^2(y_1 - 0.5)^2)\exp(-C_2^2(y_2 - 0.5)^2)
\end{equation}
where $C_1=2^{-2}C$, $C_2=3^{-2}C$ and $C=\frac{36}{13}$.
A controlled  \bkrevisesout{model}{fidelity} error is artificially introduced by considering a sequence of approximations \ah{of the quantity of interest} $q^{\miModel}:\Gamma\to\R$ for $\miModel\in\N$, \mg{such that}
\begin{equation} \label{eq:2dgpnoisy}
    {q}^{\miModel}(\y):= q(\y) + 10^{-2\miModel} x(\y)
\end{equation}
where, for each parameter $\y$, the solver error $x(\y)$ is a realisation of a \ah{standard} normal random variable $X\sim \mathcal{N}(0,1)$.
Evaluations for each \mg{fidelity} $\miModel$ \mg{have a corresponding} solve cost $\widehat{W}_{\miModel}=10^{\miModel}$.
\ah{This setup mimics a scenario in which the solver error decreases exponentially with fidelity level $\miModel$, while the  \bkrevisesout{model}{fidelity} evaluation cost increases exponentially \bkrevisesout{as $\widehat{W}_{\miModel}=10^{\miModel}$}{}.
This cost--accuracy trade off is consistent with many numerical solvers; for example, in $d$-dimensional finite element methods, \mg{uniformly} refining the mesh by a factor of $10^{\miModel}$ increases the number of degrees of freedom, and hence the computational cost, by a factor of $10^{\miModel d}$, while the approximation error decreases as $10^{-\miModel p}$ for some $p>0$ \mg{given by the accuracy of the employed discretisation}.}

MISC approximations are constructed using symmetric Leja points \cite{Piazzola2022} with the `two-step' \ah{growth} rule defined in \eqref{eq:rules}.
Errors are computed numerically with respect to a high fidelity reference \ah{solution}.
\ah{This reference $\approxRef{q}$ is obtained via} a single-fidelity \ah{sparse grid} approximation using the \bk{combination technique} \eqref{eq:combinationtechnique} with Clenshaw--Curtis points and the doubling rule of \eqref{eq:rules}, \ah{using the highest} fidelity $\approxRef{\miModel}=8$ and an isotropic Smolyak multi-index set
\begin{equation}\label{eq:isotropic-miset}
    \miset_\text{ref}=\{\vMi \st \Vert \vMi \Vert_1 \leq 2 + w \} \text{ for } w = 8.
\end{equation}
Intermediate approximations are computed at $w=0,1,2,...,7$\bk{, in \bkrevisesout{additional}{addition} to a single reference fidelity adaptive SC approximation.}
\ah{The use of reference} fidelity $\approxRef{\miModel}=8$ \ah{ensures a} reference solver error of \ah{the order of} $10^{-16}$ (i.e.\ \ah{near} machine precision) and provides a \ah{reliable benchmark for evaluating parametric convergence and model performance.}

The classical MISC response surface $\approxFinal{q}$ \mg{obtained} using \Cref{alg:misc} \ah{has} a total computational cost of $2891$ \ah{and it} is shown in \Cref{fig:2dresponsescb}.
This is clearly corrupted by the low fidelity solve error.
The \mg{corresponding} error surface is shown in \Cref{fig:2dresponsescd}. 
\begin{figure}
    \begin{subfigure}[t]{0.48\textwidth}
        \inputfrom{figures/twodgaussian/misc/}{surface_noplateau.tex}
        \caption{MISC approximation $\approxFinal{q}(\y)$.}
        \label{fig:2dresponsescb}
    \end{subfigure}
    \begin{subfigure}[t]{0.48\textwidth}
        \inputfrom{figures/twodgaussian/misc/}{errorsurface_noplateau.tex}
        \caption{MISC error surface $\approxRef{q}(\y)-\approxFinal{q}(\y)$.}
        \label{fig:2dresponsescd}
    \end{subfigure}
    \caption{MISC approximation $\approxFinal{q}$ of the 2DGP Genz test function \eqref{eq:2dgp} using \Cref{alg:misc} with fidelities $q^{\miModel}$ defined in \eqref{eq:2dgpnoisy}.}
\end{figure}
The $L_{\rho}^2(\Gamma;\R)$ approximation error is estimated as
\begin{equation}\label{eq:estimated-L2}
    \Vert q - \approxFinal{q} \Vert_{L_{\rho}^2(\Gamma;\R)} \approx \frac{1}{\lvert \Gamma \rvert} \int_{\Gamma} \lvert \approxRef{q}(\y) - \approxFinal{q}(\y) \rvert^2 \dd \y \approx \left(\frac{1}{N_{MC}}\sum_{i=1}^{N_{MC}} \lvert \approxRef{q}(\y_i) - \approxFinal{q}(\y_i) \rvert^2 \right)^{1/2}
\end{equation}
where $\rho(\y) = 1/{\lvert \Gamma \rvert}=1$ is used as the weight function in the integral.
\mg{Equation \eqref{eq:estimated-L2} is computed using $N_{MC}=10^4$ Monte Carlo samples}.
The $H_{\rho}^1(\Gamma;\R)$ norm of the approximation error is estimated similarly, using a finite difference approximation for the gradient of $\approxRef{q}$ and $\approxFinal{q}$ as implemented in the \sgmk{} \cite{Piazzola2022}.
The error plots \mg{in} \Cref{fig:ltwo2dgpmiscclassic,fig:hone2dgpmisc} show that the \bk{classical} MISC approximation constructed by \Cref{alg:misc} has not delivered an accurate approximation to the 2DGP function defined in \eqref{eq:2dgp} and the approximation quality has stagnated.
\bk{Here, and throughout the experiments, computational cost will refer to the sum of evaluation cost for all of the evaluations required by the approximation algorithm.}
The approximation is stuck adding cheap, low fidelity evaluations.
The algorithm eventually terminates, \mg{using 505 points at fidelity $\miModel=1$, resulting in the surrogate shown in \Cref{fig:2dresponsescb}}. 
\begin{figure}
    \begin{subfigure}[t]{0.48\textwidth}
        \inputfrom{figures/twodgaussian/misc/}{ltwomisc.tex}
        \caption{Estimated $L_{\rho}^{2}(\Gamma;\R)$ error.}
        \label{fig:ltwo2dgpmiscclassic}
    \end{subfigure}
    \begin{subfigure}[t]{0.48\textwidth}
        \inputfrom{figures/twodgaussian/misc/}{honemisc.tex}
        \caption{Estimated $H_{\rho}^{1}(\Gamma;\R)$ error.}
        \label{fig:hone2dgpmisc}
    \end{subfigure}
    \caption{Estimated $L_{\rho}^2(\Gamma;\R)$ and $H_{\rho}^1(\Gamma;\R)$ error for the MISC approximation $\approxFinal{q}$ of the 2DGP Genz test function \eqref{eq:2dgpnoisy} using \Cref{alg:misc} with \mg{fidelities} $q^{\miModel}$ defined in \eqref{eq:2dgpnoisy}.}
\end{figure}

The sampling profile shown in \Cref{fig:2dgpsamplingprofilemisc} shows the number of solves $\nColloc{}\rvert_{\miModel}$ for each \mg{fidelity} \bkrevise{against the cumulative computational work of the MISC algorithm. Each marker represents an iteration of the algorithm.}
\begin{figure}
    \begin{subfigure}[t]{0.48\textwidth}
        \inputfrom{figures/twodgaussian/misc/}{sampleprofile.tex}
        \caption{Number of \mg{fidelity} solves $\nColloc{}\rvert_{\miModel}$ for each \mg{fidelity} $\miModel$ \bkrevise{against cumulative computational work}.}
        \label{fig:2dgpsamplingprofilemisc}
    \end{subfigure}
    \begin{subfigure}[t]{0.48\textwidth}
        \inputfrom{figures/twodgaussian/misc/}{allenvelopes.tex}
        \caption{Spectral polynomial coefficients $\approxFinal{\FT{q}}_{\vMiPoly}\rvert_{\miModel}$ for each \mg{fidelity} $\miModel$ at final iteration.}
        \label{fig:2dgpenvelopemisc}
    \end{subfigure}
    \caption{Number of \mg{fidelity} solves against computational work and snapshot of spectral polynomial coefficients for the MISC approximation $\approxFinal{q}$ of the 2DGP Genz test function \eqref{eq:2dgpnoisy} using \Cref{alg:misc}  with \mg{fidelities} $q^{\miModel}$ defined in \eqref{eq:2dgpnoisy}.}
    \label{}
\end{figure}
Many more grid points \mgrevise{$\nColloc{\miset} \vert_{\miModel}$} are sampled for fidelity $\miModel=1$, \mgrevise{whilst} only one grid point \mgrevise{is sampled} for fidelity $\miModel=3$.
By definition in \eqref{eq:2dgpnoisy}, it is clear that the \mg{fidelity} $\miModel=1$ is extremely erroneous.
\Cref{fig:2dgpenvelopemisc} shows the spectral polynomial coefficients $\{\approxFinal{\FT{q}}_{\vMiPoly}\rvert_{\miModel} \}_{\vMiPoly\in\misetPoly(\miset)}$ for each \mg{fidelity} in the computed MISC approximation $\approxFinal{q}$.
Notice both the \mg{fidelities} $\miModel=1$ and $\miModel=2$ have notable \emph{spectral plateaus} in the polynomial coefficients.
For fidelity $\miModel=1$, it is clear that the solver error is spoiling the computed spectral polynomial coefficients, and consequently also \mg{leading to an} approximation \mg{of poor} quality.
The plateau is lowered by using a more accurate solver, \mg{with} $\miModel=2$.
There is only one solve on the \mg{fidelity} $\miModel=3$ so there is only the constant polynomial $\vMiPoly=\mistyle{0}$ plotted in \Cref{fig:2dgpenvelopemisc}.


\subsubsection{Plateau detection}
Plateau detection via \Cref{alg:plateau} is applied to the fidelities in the classical MISC approximation with \ah{$\nBurnIn=\nBurnOut=2$}. and $k_{min}=3$.
The envelopes for the three active \mg{fidelities} $\miModel=1,2,3$ are shown in \Cref{fig:2dgpenv}.
For $\miModel=1,2$, these capture the exponential decrease in coefficient magnitude and the plateau in the coeffic\bkrevise{i}ent magnitude due to solver noise.
The piecewise log-linear models are shown in \Cref{fig:2dgpplat}.
The effect of the burn in $\nBurnIn$ and burn out $\nBurnOut$ are observed.
The first models do not start at zero, but at polynomial total degree equal to two.
Similarly, the second models finish two units before the end of the envelopes.
The slopes of the second piecewise log--linear models are $-8.0\times10^{-4}$ and $-1.4\times10^{-2}$ for \mg{fidelities} $\miModel=1,2$ respectively.
In both cases a plateau is successfully identified by \Cref{alg:plateau}.
The computed plateau levels are $\error{plateau}=4.5\times10^{-3}$ and $8.4\times10^{-5}$ respectively.
This matches intuition: the solver error defined in \eqref{eq:2dgpnoisy} for \mg{fidelity} $\miModel=1$ has standard deviation $10^{-2}$ and for \mg{fidelity} $\miModel=2$ has standard deviation $10^{-4}$.
\bk{For \mg{fidelity} $\miModel=3$ the envelope \bkrevise{in \Cref{fig:2dgpplat}} is not visible --- it is a single point, \bkrevise{corresponding to the single polynomial coefficient for $\miModel=3$}!}

\begin{figure}
    \begin{subfigure}[t]{0.48\textwidth}
        \inputfrom{figures/twodgaussian/misc/}{plateaudetection_env.tex}
        \caption{Spectral polynomial coefficients $\approxFinal{\FT{q}}_{\vMiPoly}$ and envelopes $\envelope \rvert_{\miModel}$ for each \mg{fidelity} $\miModel$
            at final iteration. \bkrevise{The plot for fidelity $\alpha=3$ consists of just of a single marker at $\Vert \vMiPoly \Vert_1 = 0$.}}
        \label{fig:2dgpenv}
    \end{subfigure}
    \begin{subfigure}[t]{0.48\textwidth}
        \inputfrom{figures/twodgaussian/misc/}{plateaudetection_plateau.tex}
        \caption{Envelopes $\envelope \rvert_{\miModel}$ (dashed) and piecewise log-linear models (solid) for each \mg{fidelity} $\miModel$ at final iteration.
            \bkrevise{Note that although fidelity $\miModel=3$ is included in the legend, the envelope consists of a single point at $\Vert \vMiPoly \Vert_1 = 0$
                so it is not visible.}} \label{fig:2dgpplat}
    \end{subfigure}
    \caption{Spectral polynomial coefficients, spectral envelopes and piecewise log-linear models for the MISC approximation $\approxFinal{q}$ of the 2DGP Genz test function \eqref{eq:2dgpnoisy} using \Cref{alg:misc} with \mg{fidelities} $q^{\miModel}$ defined in \eqref{eq:2dgpnoisy}.}
    \label{fig:2dgp-plateaus}
\end{figure}

\subsubsection{\plateaumisc{}}
\mg{The} \plateaumisc{} \mg{strategy in} \Cref{alg:pmisc}, is now applied: comparing the resulting response surface (\Cref{fig:2dresponses-pmisc}) to the MISC \mg{surrogate} in \Cref{fig:2dresponsescb} it is clear \mg{that} \Cref{alg:pmisc} offers an improved response surface --- it is far less noisy.
\mg{The corresponding error surface  $\approxRef{q}(\y)-\approxFinal{q}(\y)$ provided by \plateaumisc{} is reported in \Cref{fig:2dresponses-pmisc-error}}.
\begin{figure}
    \begin{subfigure}[t]{0.48\textwidth}
        \inputfrom{figures/twodgaussian/misc/}{surface.tex}
        \caption{\plateaumisc{} approximation $\approxFinal{q}(\y)$.}
        \label{fig:2dresponses-pmisc}
    \end{subfigure}
    \begin{subfigure}[t]{0.48\textwidth}
        \inputfrom{figures/twodgaussian/misc/}{errorsurface.tex}
        \caption{\plateaumisc{} error surface $\approxRef{q}(\y)-\approxFinal{q}(\y)$.}
        \label{fig:2dresponses-pmisc-error}
    \end{subfigure}
    \caption{\plateaumisc{} approximation $\approxFinal{q}$ of the 2DGP Genz test function \eqref{eq:2dgp} using \Cref{alg:pmisc} with fidelities $q^{\miModel}$ defined in \eqref{eq:2dgpnoisy}.}
    \label{fig:2dgp-responses-pmisc}
\end{figure}

To make a fair comparison, the approximation error is plotted against computational cost for both approximations.
The estimated $L^2_\rho(\Gamma;\R)$ and $H^1_\rho(\Gamma;\R)$ errors are shown in \Cref{fig:ltwo2dgp-pmisc,fig:hone2dgp-pmisc}.
\begin{figure}
    \begin{subfigure}[t]{0.48\textwidth}
        \inputfrom{figures/twodgaussian/misc/}{ltwo-pmisc.tex}
        \caption{Estimated $L_{\rho}^{2}(\Gamma;\R)$ error}
        \label{fig:ltwo2dgp-pmisc}
    \end{subfigure}
    \begin{subfigure}[t]{0.48\textwidth}
        \inputfrom{figures/twodgaussian/misc/}{hone-pmisc.tex}
        \caption{Estimated $H_{\rho}^{1}(\Gamma;\R)$ error}
        \label{fig:hone2dgp-pmisc}
    \end{subfigure}
    \caption{Estimated $L_{\rho}^2(\Gamma;\R)$ and $H_{\rho}^1(\Gamma;\R)$ approximation error for the \plateaumisc{} approximation $\approxFinal{q}$ of the 2DGP Genz test function \eqref{eq:2dgpnoisy} using \Cref{alg:pmisc}  with \mg{fidelities} $q^{\miModel}$ defined in \eqref{eq:2dgpnoisy}. Vertical arrows denote the times at which a fidelity is marked as saturated.}
\end{figure}
The \plateaumisc{} approximation error decreases as we increase the computational budget since it can successfully transition to higher available fidelities, whereas the classical MISC approximation stagnates at a large error and cannot proceed further as it requests too many points from the noisy low fidelities.
Vertical arrows in \Cref{fig:ltwo2dgp-pmisc,fig:hone2dgp-pmisc} (and later figures) indicate points at which a fidelity is marked as saturated in the \plateaumisc{} approximation.
These are generally aligned with a decrease in the approximation error, however the approximation is not returned at every iteration of the algorithm (output \bk{iterations} are controlled by computational budget) so there is a discrepancy at some points.
The benefit of the multi-fidelity approach is also shown.
For comparable computational cost, a more accurate approximation is constructed by \plateaumisc{} than the single high fidelity reference  \bkrevisesout{model}{} or a single reference fidelity adaptive stochastic collocation approximation.

The final multi-index sets are shown in \Cref{fig:miset-2dgp}: it is clear \bk{that} MISC oversamples the low \bkrevisesout{fidelity models}{fidelities} and does not access the more accurate \mg{fidelities} with higher $\miModel$.
The effect of backfilling \mg{in \plateaumisc{}} is also clearly observed: \bk{the lower fidelity levels (smaller $\miModel$) in the \plateaumisc{} multi-index set are enriched to allow the higher \bkrevisesout{fidelity models}{fidelities} to access higher polynomial degree parametric approximation spaces}.
For problems with stronger correlations in the errors \bk{between fidelities}, the multi-index set will be more anisotropic with respect to fidelity \bk{since more information can be captured by the low fidelities}.
\mg{This will be observed later in \Cref{sec:parabolic}}.
\begin{figure}
    \begin{subfigure}{0.49\textwidth}
        \inputfrom{figures/twodgaussian/misc/}{miset_np.tex}
        \caption{Multi-index set $\miset$ for final MISC approximation}
    \end{subfigure}
    \begin{subfigure}{0.49\textwidth}
        \inputfrom{figures/twodgaussian/misc/}{miset.tex}
        \caption{Multi-index set $\miset$ for final \plateaumisc{} approximation}
    \end{subfigure}
    \caption{Multi-index sets $\miset$ for the final MISC and \plateaumisc{} approximations of the 2DGP Genz test function \eqref{eq:2dgpnoisy} using \Cref{alg:misc,alg:pmisc} with \mg{fidelities} $q^{\miModel}$ defined in \eqref{eq:2dgpnoisy}.}
    \label{fig:miset-2dgp}
\end{figure}

From the fidelity work allocations reported in \Cref{fig:2dgpsamplingprofile-pmisc}, we see that computational cost is saved in \plateaumisc{} as the high \bkrevisesout{fidelity models}{fidelities} are sequentially activated, and only when the parametric approximation is sufficiently accurate for it to be worthwhile.
The effect of backfilling is seen as the number of solves will continue to increase for all fidelities, \bk{despite some being marked as saturated and blocked from explicit refinement}.
Both the exponential decrease in the spectral coefficient magnitudes and the plateau effect from the solver error are clearly observed in the spectral polynomial coefficients in \Cref{fig:2dgpenvelope-pmisc}.
\begin{figure}
    \begin{subfigure}[t]{0.48\textwidth}
        \inputfrom{figures/twodgaussian/misc/}{sampleprofile-pmisc.tex}
        \caption{Number of \mg{fidelity} solves $\nColloc{}\rvert_{\miModel}$ for each \mg{fidelity} $\miModel$ \bkrevise{against cumulative computational work}.}
        \label{fig:2dgpsamplingprofile-pmisc}
    \end{subfigure}
    \begin{subfigure}[t]{0.48\textwidth}
        \inputfrom{figures/twodgaussian/misc/}{plateaudetection_env_pmisc.tex}
        \caption{Spectral envelope $\envelope\rvert_{\miModel}$ for each \mg{fidelity} $\miModel$ at \bkrevise{plateau detection}
            \bkrevise{(which happens at different iterations for  different fidelities)}.}
        \label{fig:2dgpenvelope-pmisc}
    \end{subfigure}
    \caption{Number of \mg{fidelity} solves against computational work and snapshot of spectral polynomial coefficients for the \plateaumisc{} approximation $\approxFinal{q}$ of the 2DGP Genz test function \eqref{eq:2dgpnoisy} using \Cref{alg:pmisc} with \mg{fidelities} $q^{\miModel}$ defined in \eqref{eq:2dgpnoisy}.}
\end{figure}

Finally, we consider approximating the probability density function (PDF) $\rho_{q}:\Gamma\to\R$ of the quantity of interest for uniformly distributed input parameters, that is estimating $\rho_{q}$ such that
\begin{equation}
    \mathrm{Prob}(q \in A) = \int_{A} \rho_{q}(\y) \dd \y.
\end{equation}
The PDF is estimated using \mg{$N_{MC}=10^4$} Monte Carlo samples from the surrogate models and the \matlab{} \texttt{ksdensity} function.
\begin{figure}
    \begin{subfigure}[t]{0.49\textwidth}
        \inputfrom{figures/twodgaussian/misc/}{pdf.tex}
        \caption{Approximate probability density functions for the final approximations}
        \label{fig:2dgp-pdf}
    \end{subfigure}
    \begin{subfigure}[t]{0.49\textwidth}
        \inputfrom{figures/twodgaussian/misc/}{ks.tex}
        \caption{Evolution of two sample Kolmogorov--Smirnov (KS2) test statistic with computational cost}
        \label{fig:2dgp-ks}
    \end{subfigure}
    \caption{Final approximate probability density functions and evolution of two sample Kolmogorov--Smirnov test statistic for the 2DGP Genz test function \eqref{eq:2dgpnoisy} using \Cref{alg:pmisc} with \mg{fidelities} $q^{\miModel}$ defined in \eqref{eq:2dgpnoisy}.}
\end{figure}
\Cref{fig:2dgp-pdf} \mg{shows that} the \plateaumisc{} PDF approximation closely matches the reference PDF approximation whereas the classical MISC PDF approximation is \mg{significantly less accurate}.
This is \mg{also} supported by the two-sample Kolmogorov--Smirnov test statistic, as shown in \Cref{fig:2dgp-ks}.
This is produced using $N_{MC}=10^6$ samples from the surrogates and the reference approximation.
Even at relatively low computation cost, the \plateaumisc{} approximation gives a good approximation \mg{of} the PDF, whilst the classical MISC approximation shows no convergence with computational cost at all.
The plateau in the KS2 test statistic is hypothesised to be a consequence of discretisation due to the finite sample sets.

\section{Numerical examples} \label{sec:numerics}
\ah{To illustrate the practical relevance of the proposed \plateaumisc{} method we consider two representative test cases: an \mg{unsteady} advection--diffusion problem \mg{with uncertain coefficients}, and a \mg{parametric} turbulent incompressible Navier--Stokes flow.
    The first example is a parametric variant of the well-known \emph{double glazing problem}, originally introduced in \cite{Elman2014,Elman2014b} and based on the benchmark problems from the International Association for Hydraulic Research \cite{Smith1982,Morton1996}.}
\mgrevise{This problem features a two-dimensional domain in space, a time dependence, and  two uncertain parameters, leading to a high-dimensional problem of dimension 5.}
The second example involves a parametric \ah{version of the \emph{Two Dimensional Wall Mounted Hump (2DWMH)} problem,} proposed by NASA for \ah{turbulence model} validation \cite{NASA2DWMH,Greenblatt2004}.
\ah{This \mgrevise{highly non-linear problem has been extended in \cite{Tsiolakis2022} by introducing}} a parametric suction jet \mg{to control} the reattachment point \ah{of the separated flow downstream of} the `hump'.
\mgrevise{In our work, the Reynolds number is also introduced as a second parameter, leading to a problem set in a high-dimensional space of dimension 4.}

In both cases, the fidelity is controlled by a single \bkrevisesout{dimension}{index} ($\nModel=1$)
\bkrevisesout{while the parametric domain is two dimensional ($\nY=2)$.
Although sparse grid methods are typically applied to high-dimensional problems, the focus here is managing the impact of the solver error,
which is equally significant in these lower-dimensional yet practically relevant settings.
}{\MGnew{associated with} the accuracy of the solver employed for the computation of each instance of the parametric problem, }
\mgrevise{while the parametric dimension is set to $d=2$. It is worth noticing that MISC has already been successfully validated for problems with parametric dimension larger than $d=2$, see e.g. \cite{HajiAli2016a,HajiAli2016b,Jakeman2019}, where experimental tests for up to tens, hundreds, and thousands of parameters were respectively presented.
The scope of this work is instead to investigate the impact of the solver error on the accuracy of the surrogate model of the QoI and the development of an automatic strategy to make  MISC robust in this framework.
The presented numerical experiments thus focus on exploring different sources of such error and their interaction with nonlinearity, rather than on high-dimensional problems.
We point out that the plateau detection algorithm is essentially dimension-independent in terms of cost, so its execution as parameter dimension $d$ increases does not pose any additional challenge compared to the traditional versions of MISC.}


\subsection{Uncertainty quantification for unsteady advection--diffusion problem with uncertain coefficients} \label{sec:parabolic}
Consider the spatial domain $\spatialdomain:=(-1,1)^2$, finite time $T>0$ and parameter domain $\Gamma:=[0,1]^2$.
\mg{For each $\y\in\Gamma$, we seek the solution to the following advection--diffusion problem:}
\begin{equation}\label{eq:parabolic}
    \begin{cases}
        \displaystyle\frac{\partial u}{\partial t}(\x,t;\y) - \nabla \cdot \big(\diffFn(\y) \nabla u(\x,t;\y)\big) + \wind(\x,\y) \cdot \nabla u(\x,t;\y) = 0, & (\x,t) \in \spatialdomain \times (0,T), \\
        u(\x,t;\y) = (1-\exp(-t/0.1))(1+0.1\sin(2\pi t)),                                                                                                      & (\x,t) \in \spatialbdry \times (0,T),   \\
        u(\x,0;\y) = 0,                                                                                                                                        & \x \in \overline{\spatialdomain},
    \end{cases}
\end{equation}
where the operator $\nabla$ acts on the spatial variable $\x$ and \ah{the diffusion coefficient $\diffFn(\y)$ and advection field $\wind(\x,\y)$ are affine functions of the parameters $\y\in\Gamma$, \mg{namely},}
\begin{equation}\label{eq:windanddiff}
    \diffFn(\y) :={} 0.1\diffFn_0 + 1.8\diffFn_0 y_1, \quad
    \wind(\x,\y) :={} \wind_0(\x) + \sum_{i=1}^{4} \wind_i(\x) (2y_2-1),
\end{equation}
with
\begin{equation}
    \begin{aligned}
        \wind_0(\x) & :={} [2x_2(1-x_1^2), -2 x_2 (1-x_1^2)], \quad
        \wind_i(\x) & :={}
        \begin{cases}
            \wind_0(2(\x-[c_i,d_i])) & \x \in\spatialdomain_i, \\
            [0,0],                   & \text{otherwise},
        \end{cases}
    \end{aligned}
\end{equation}
for $\vec{c}:=[-0.5,-0.5,0.5,0.5]$, $\vec{d}:=[-0.5,0.5,-0.5,0.5]$, $\spatialdomain_i:=[c_i-0.5, d_i+0.5]\times[c_i-0.5, d_i+0.5]$.
Three realisations of the advection field $\wind(\x,\y)$ are given in \Cref{fig:windrealisations}.
\begin{figure}[t]
    \hfill\begin{subfigure}[t]{0.33\textwidth}
        \inputfrom{figures/dg/}{advectionfield_-1.tex}
        \caption{$\wind\big(\x,[-1,-1]\big)$}
    \end{subfigure}\hfill
    \begin{subfigure}[t]{0.33\textwidth}
        \pgfmathdeclarefunction{yval}{0}{%
            \pgfmathparse{0.0}
        }
        \inputfrom{figures/dg/}{advectionfield_0.tex}
        \caption{$\wind\big(\x,[0,0]\big)$}
    \end{subfigure}\hfill
    \begin{subfigure}[t]{0.33\textwidth}
        \pgfmathdeclarefunction{yval}{0}{%
            \pgfmathparse{1.0}
        }
        \inputfrom{figures/dg/}{advectionfield_1.tex}
        \caption{$\wind\big(\x,[1,1]\big)$}
    \end{subfigure}\hfill
    \caption{Realisations of the advection field $\wind\big(\x,\y\big)$ defined in \eqref{eq:windanddiff} \mg{for} three parameter realisations $\y=[-1,-1],[0,0]$ and $[1,1]$.}
    \label{fig:windrealisations}
\end{figure}

Equation \eqref{eq:parabolic} has a parametric diffusion coefficient and a parametric perturbation to the average advection field $\wind_0$.
The solution does not reach a steady--state due to the time dependence in the Dirichlet boundary condition.
Observe that for $\diffFn_0>0$ the problem is well posed because the diffusion coefficient is bounded $0.09\diffFn_0<\diffFn(\y)<2\diffFn_0$ and the advection field is finite for all parameters $\y\in\Gamma$, but never advection dominated.
Further discussion on the well-posedness of this problem and verification of the analyticity of the solution $u(\x,t;\y)$ with respect to the parameters $\y$ can be found in \cite{Kent2023b}.

The quantity of interest is $q(\y):=u(\x_*,t_*;\y)$ for the space-time coordinate $(\x_*,t_*):=([0.5,-0.5],10)$.
\ah{Note that the chosen QoI is a pointwise quantity rather than an average or integrated quantity.
    This introduces additional difficulty as pointwise QoIs are typically more sensitive to local solution features and numerical errors, making accurate approximation more challenging.}

\subsubsection{Numerical approximation of PDE solution}
Pointwise (in $\Gamma$) approximations of the QoI $q(\y)$ are computed by using a $\meshP_1$ \mg{continuous Galerkin} finite element (FE) approximation of \eqref{eq:parabolic} (piecewise linear polynomials on a structured uniform triangular mesh) and solving the resulting system of ODEs using \trab{} adaptive timestepping with local error control \cite{Iserles2008,Gresho2008}. 
The solver accuracy (fidelity) is controlled by the local timestepping error tolerance $\letol>0$ whilst the mesh is fixed with $400\times400$ uniformly spaced \mg{points}.
The \ah{maximum} mesh P\'eclet number $\Peclet{}^{h}$ is \ah{obtained} for the diffusion coefficient $\diffFn\big([0,y_2]\big)=0.1 \diffFn_0$ and maximum advection field magnitude with $y_2=1$.
This can be bound by \ah{$\Peclet{}^{h}_{max} = {(h \lvert \wind_{max}\rvert)}/{(2 {\diffFn_{min}})} = {(4 \cdot 2 \cdot 400^{-1})}/{(2 \cdot 10^{-1}\cdot 10^{-1})} = 1$} \mg{and it has been numerically verified that the mesh provides stable solutions even in the worst case scenario}.
Cost is proportional to the number of timesteps used to reach $t_*=10$ which is approximately inversely proportional to the size of timesteps.
For \trab{} the timesteps are approximately proportional to $\letol^{1/3}$ so the solver cost is modelled as
\begin{equation}\label{eq:parabolicwork}
    \widehat{W}_{\letol} = \letol^{-1/3}.
\end{equation}

There are two main sources of solver error: spatial and temporal.
\ah{The spatial error in a semi-discrete solution using a conforming finite element approximation does not change the regularity of the computed approximation \cite{Nobile2009,Adcock2022}.
    Conversely, the timestepping error does \mg{affect} the regularity of the solution with respect to the parameter $\y$, and \mg{it} is the main source of the spectral coefficient plateau in our computed approximation.}


\subsubsection{Multi-fidelity MISC approximation}
The local timestepping error \mg{(LTE)} tolerance $\letol$ is used to define the \mg{different fidelities for the problem such that}
\begin{equation}\label{eq:models-parabolic}
    \miModel \mapsto \{\approxS{q}:\Gamma\to\R \text{ computed with LTE tolerance } \letol=10^{-\miModel}\}
\end{equation}
with associated \bkrevisesout{model}{} evaluation cost $\widehat{W}_{\miModel}=10^{\miModel/3}$ derived from \eqref{eq:parabolicwork}.
Symmetric Leja points and the two-step rule from \eqref{eq:rules} are used to form the MISC approximation.

As in \Cref{example:2dgp}, a reference approximation is constructed using an isotropic Smolyak sparse grid with multi-index set defined by \eqref{eq:isotropic-miset} with $w=5$ and Clenshaw--Curtis points with the doubling rule from \eqref{eq:rules}.
The reference approximation uses the  \bkrevisesout{model}{fidelity} defined by $\approxRef{\miModel}=7$.
Three realisations of $\approxRef{u}(\cdot,t_*;\y)$ for parameters $\y=[-1,-1],[0,0]$ and $[1,1]$ are shown in \Cref{fig:dgrealisations}.
\begin{figure}[t]
    \hfill\begin{subfigure}[t]{0.33\textwidth}
        \inputfrom{figures/dg/}{sample-1.tex}
        \caption{$\approxRef{u}(\x,10;[-1,-1])$}
    \end{subfigure}\hfill
    \begin{subfigure}[t]{0.33\textwidth}
        \inputfrom{figures/dg/}{sample0.tex}
        \caption{$\approxRef{u}(\x,10;[0,0])$}
    \end{subfigure}\hfill
    \begin{subfigure}[t]{0.33\textwidth}
        \inputfrom{figures/dg/}{sample1.tex}
        \caption{$\approxRef{u}(\x,10;[1,1])$}
    \end{subfigure}\hfill
    \caption{Realisations of the reference approximation $\approxRef{u}(\x,10;\y)$ to the parametric parabolic PDE problem \eqref{eq:parabolic} at three parameter realisations $\y=[-1,-1],[0,0]$ and $[1,1]$. The marker $\times$ indicates the spatial location $\x_*=(0.5,-0.5)$ where the QoI is evaluated.}
    \label{fig:dgrealisations}
\end{figure}
The effects of changing the diffusion coefficient $\diffFn(\y)$ and the advection field $\wind(\cdot,\y)$ are clearly illustrated.
For example, as the diffusion coefficient $\diffFn(\y)$ increases, the diffusive boundary layer on the right hand wall becomes wider.
Similarly, as the advection field $\wind(\cdot,\y)$ increases in strength, the fluid is transported more strongly into the domain.
Histograms of the approximation error \ah{of the QoI} for $\miModel=1,2,3$ are shown in \Cref{fig:hist-dg}.
It can be seen \mg{that} as fidelity \ah{level} $\miModel$ increases, the approximation error decreases but the distribution is not necessarily Gaussian.
\begin{figure}
    \begin{subfigure}{0.33\textwidth}
        \inputfrom{figures/dg/}{histogram_1.tex}
        \caption{\mg{Fidelity} $\miModel=1$}
    \end{subfigure}
    \begin{subfigure}{0.33\textwidth}
        \inputfrom{figures/dg/}{histogram_2.tex}
        \caption{\mg{Fidelity} $\miModel=2$}
    \end{subfigure}
    \begin{subfigure}{0.33\textwidth}
        \inputfrom{figures/dg/}{histogram_3.tex}
        \caption{\mg{Fidelity} $\miModel=3$}
    \end{subfigure}
    \caption{Histograms of approximation error $q^\text{ref} - q^{\miModel}$ from $n=1000$ Latin Hypercube samples for the parametric parabolic PDE \eqref{eq:parabolic} with \mg{fidelities} $q^{\miModel}$ defined in \eqref{eq:models-parabolic}.}
    \label{fig:hist-dg}
\end{figure}

Computed response surfaces using the \plateaumisc{} and MISC algorithms are shown in \Cref{fig:miscdgsurg} alongside pointwise error surfaces with respect to the reference approximation.
\begin{figure}
    \begin{subfigure}[t]{0.48\textwidth}
        \inputfrom{figures/dg/misc/}{surface_noplateau.tex}
        \caption{\textsc{MISC} response surface $\approxFinal{q}(\y)$.}
        \label{fig:dgclassicmiscsurf}
    \end{subfigure}
    \begin{subfigure}[t]{0.48\textwidth}
        \inputfrom{figures/dg/misc/}{surface.tex}
        \caption{\plateaumisc{} response surface $\approxFinal{q}(\y)$.}
        \label{fig:dgmiscsurf}
    \end{subfigure}
    \begin{subfigure}[t]{0.48\textwidth}
        \inputfrom{figures/dg/misc/}{errorsurface_noplateau.tex}
        \caption{\textsc{MISC} error surface $\approxRef{q}(\y)-\approxFinal{q}(\y)$.}
        \label{fig:errorclassicmiscdg}
    \end{subfigure}
    \begin{subfigure}[t]{0.48\textwidth}
        \inputfrom{figures/dg/misc/}{errorsurface.tex}
        \caption{\plateaumisc{} error surface $\approxRef{q}(\y)-\approxFinal{q}(\y)$.}
        \label{fig:errormiscdg}
    \end{subfigure}
    \caption{Approximated response surfaces $\approxFinal{q}(\y)$ for the parametric parabolic PDE \eqref{eq:parabolic} using (a) \Cref{alg:misc} and (b) \Cref{alg:pmisc} with \mg{fidelities} $q^{\miModel}$ defined in \eqref{eq:models-parabolic}.
        \Cref{fig:errormiscdg,fig:errorclassicmiscdg} illustrate the \mg{corresponding} approximation error $\approxRef{q}(\y)-\approxFinal{q}(\y)$.}
    \label{fig:miscdgsurg}
\end{figure}
\begin{figure}
    \begin{subfigure}[t]{0.48\textwidth}
        \inputfrom{figures/dg/misc/}{ltwo.tex}
        \caption{Estimated $L^{2}(\Gamma;\R)$ error}
    \end{subfigure}
    \begin{subfigure}[t]{0.48\textwidth}
        \inputfrom{figures/dg/misc/}{hone.tex}
        \caption{Estimated $H^{1}(\Gamma;\R)$ error}
    \end{subfigure}
    \caption{Estimated $L_{\rho}^2(\Gamma;\R)$ and $H_{\rho}^1(\Gamma;\R)$ approximation error for the \plateaumisc{} approximation $\approxFinal{q}$ of parametric parabolic PDE \eqref{eq:parabolic} using \Cref{alg:pmisc} with \mg{fidelities} $q^{\miModel}$ defined in \eqref{eq:models-parabolic}. Vertical arrows denote the times at which a fidelity is marked as saturated.}
    \label{fig:dgerror}
\end{figure}
The pointwise error in the \plateaumisc{} approximation is generally more than one order of magnitude smaller than in the MISC approximation.

The benefits of \plateaumisc{} are shown more clearly in the error plots in \Cref{fig:dgerror}.
Note that the MISC and \plateaumisc{} approximations are not identical even for low computational cost: this is a consequence of constructing the output \plateaumisc{} approximation using only the multi-index set $\miset$, rather than the larger multi-index set \bk{$\miset\cup\reducedmargin{\miset}$}.
The surfaces in \Cref{fig:miscdgsurg} correspond to the last markers in \Cref{fig:dgerror}.
The estimated $L_{\rho}^2(\Gamma;\R)$ error is computed using \eqref{eq:estimated-L2} with $N_{MC}=10^4$ samples.
For the adaptive MISC algorithm, the $L_{\rho}^2(\Gamma;\R)$ error stagnates at approximately $10^{-2}$ whilst the \plateaumisc{} approximation error decreases to $3\times10^{-4}$.
Similarly, the $H_{\rho}^1(\Gamma;\R)$ error for the adaptive MISC algorithm stagnates at approximately $10^{0}$ whilst the \plateaumisc{} approximation error decreases to $10^{-2}$.
It is clear that the \plateaumisc{} algorithm gives a better approximation than the standard adaptive MISC approximation.
The approximation errors are also compared to the sequential construction of the reference approximation with multi-index sets \ah{$\{\vMiParam \in\N^2 \st  \Vert \vMiParam \Vert_1 \leq 2+w \}$ for $w=0,1,2,3,4$ and fidelity level $\approxRef{\miModel}=7$}. 
The \plateaumisc{} approximation attains the same rate of convergence as the reference approximation, whilst \mg{achieving} comparable approximation error for a reduced computational cost.
This is because \plateaumisc{} can balance the solver error and parametric error giving a more efficiently constructed approximation. 
\bk{Similarly, using \plateaumisc{} we see an improvement over the adaptive single reference fidelity approximation.}

\ah{The reduced computational cost is due to the balance between the solver error and parametric error.}
Initially, low \bkrevisesout{fidelity models}{fidelities} can be used as the dominant error is the \ah{parametric} \ah{interpolation error.
    In \Cref{fig:dgsamplingprofile-pmisc} we see that initially only fidelity $\miModel=1$ is sampled.}
As the total approximation error decreases, the dominant error becomes the solver error \ah{(i.e.\ time discretisation error)}.
\plateaumisc{} \mg{detects} this \mg{change}, allowing higher \bkrevisesout{fidelity models}{fidelities} to be accessed and reducing the overall approximation error further.
\ah{The sequential activation of higher fidelities is shown in \Cref{fig:dgsamplingprofile-pmisc}.}
This is not the case for classic MISC \bk{as it is harder for higher fidelities to be activated}, resulting in the stagnation of the error observed in \Cref{fig:dgerror}.

\bkrevise{The final detected set of saturated fidelities is $S=\{1,2,3,4\}$.
    The spectral plateau effect is \bkrevise{observed} for fidelities $\miModel=1,2,3,4$ in the spectral envelopes in \Cref{fig:dgenvelope-pmisc}.
    The benefit of including the burn out parameter is also clear: in detecting the plateau for $\miModel=3$ it is clearly advantageous to drop the final two coefficients in the envelope.
    The spectral envelopes for fidelities \MGnew{$\miModel=5$ and $\miModel=6$} continue to decay, and the adaptive algorithm would continue to explore these fidelities.}
\begin{figure}
    \begin{subfigure}[t]{0.48\textwidth}
        \inputfrom{figures/dg/}{sampleprofile-pmisc.tex}
        \caption{Number of solves $\nColloc{}\rvert_{\miModel}$ for each \mg{fidelity} $\miModel$ at each iteration.}
        \label{fig:dgsamplingprofile-pmisc}
    \end{subfigure}
    \begin{subfigure}[t]{0.48\textwidth}
        \inputfrom{figures/dg/}{plateaudetection_env_pmisc.tex}
        \caption{\bkrevise{Solid lines denote }spectral envelopes $\envelope\rvert_{\miModel}$ for each \mg{fidelity} $\miModel$ at \bkrevise{plateau detection}
            \bkrevise{(which happens at different iterations for different fidelities)}.
            \bkrevise{Dashed lines denote spectral envelopes for $\miModel=5,6$ at the final iteration, with no spectral plateau detected.}}
        \label{fig:dgenvelope-pmisc}
    \end{subfigure}
    \caption{Number of \mg{fidelity} solves against computational work and snapshot of spectral polynomial coefficients for the \plateaumisc{} approximation $\approxFinal{q}$ of the parametric parabolic PDE \eqref{eq:parabolic} using \Cref{alg:pmisc} with \mg{fidelities} $q^{\miModel}$ defined in \eqref{eq:models-parabolic}.}
    \label{fig:dgsamplesenvelopes}
\end{figure}

\mgrevise{
\begin{remark}
The spectral envelopes in \Cref{fig:dgenvelope-pmisc} are represented at the iteration in which the plateau is first detected.
After stagnation of a fidelity is detected, the algorithm can still add further collocation points to said fidelity, due to the back filling procedure using the multi-index set defined in \eqref{eq:backfill}.
After adding additional collocation points, the corresponding spectral envelope may no longer show a plateau, but having been previously detected by \plateaumisc{} it remains marked as saturated.
Finally, note that no spectral plateau is detected for fidelities \MGnew{$\miModel=5$ and $\miModel=6$. Hence,} the corresponding spectral envelopes (denoted by a dashed line in \Cref{fig:dgenvelope-pmisc}) continue to decrease, with the associated solver noise not affecting the quality of the constructed surrogate model.
\end{remark}
}

\ah{The multi-index sets $\miset$ at the final algorithm iterations are shown in \Cref{fig:miset-dg-both}.
    \Cref{fig:miset-dg-noplateau} indicates that, in the classic MISC algorithm, $\miModel=1,2$ are erroneously oversampled.
    The \plateaumisc{} algorithm prevents this, with the resulting multi-index set shown in \Cref{fig:miset-dg}.}
The `sloped' shape in the dimension $\miParam_2$ represents the lower fidelity information that is incorporated into the final approximation\bk{: higher degree polynomial terms are used in fidelities $\miModel=1,2,3$ than $\miModel=4,5,6$.}
\begin{figure}
    \begin{subfigure}{0.49\textwidth}
        \inputfrom{figures/dg/misc/}{miset_np.tex}
        \caption{Multi-index set $\miset$ for final MISC approximation}
        \label{fig:miset-dg-noplateau}
    \end{subfigure}
    \begin{subfigure}{0.49\textwidth}
        \inputfrom{figures/dg/misc/}{miset.tex}
        \caption{Multi-index set $\miset$ for final \plateaumisc{} approximation}
        \label{fig:miset-dg}
    \end{subfigure}
    \caption{Multi-index sets $\miset$ for the final MISC and \plateaumisc{} approximations of the parametric parabolic PDE \eqref{eq:parabolic} using \Cref{alg:misc,alg:pmisc} with \mg{fidelities} $q^{\miModel}$ defined in \eqref{eq:models-parabolic}.}
    \label{fig:miset-dg-both}
\end{figure}

The PDFs are again approximated using $N_{MC}=10^4$ samples with the \matlab{} \texttt{ksdensity} function and the KS2 test statistic is computed with $N_{MC}=10^6$ samples.
In \Cref{fig:dg-pdfs}, the distinction between MISC and \plateaumisc{} \bk{PDFs} is less clear than in the analytic example of \Cref{example:2dgp}, but \plateaumisc{} still outperforms classic MISC \mg{in terms of accuracy, efficiency and robustness.}
\begin{figure}
    \begin{subfigure}[t]{0.49\textwidth}
        \inputfrom{figures/dg/misc/}{pdf.tex}
        \caption{Approximate probability density functions for the final approximations}
    \end{subfigure}
    \begin{subfigure}[t]{0.49\textwidth}
        \inputfrom{figures/dg/misc/}{ks.tex}
        \caption{Evolution of two sample Kolmogorov--Smirnov (KS2) test statistic with computational cost}
    \end{subfigure}
    \caption{Final approximate probability density functions and evolution of two sample Kolmogorov--Smirnov test statistic for the parametric parabolic PDE \eqref{eq:parabolic} using \Cref{alg:pmisc} with  \bkrevisesout{models}{fidelities} $q^{\miModel}$ defined in \eqref{eq:models-parabolic}.}
    \label{fig:dg-pdfs}
\end{figure}

\subsection{Parametric flow control in turbulent incompressible Navier--Stokes problem}
\label{example:2dwmh}
This example considers the NASA two dimensional wall mounted hump (2DWMH) benchmark test \ah{case} \cite{NASA2DWMH,Greenblatt2004,Naughton2004}, \ah{which simulates a turbulent incompressible flow over a Glauert--Goldschmied-type body} \mg{controlled by a parameter $\y\in\Gamma=[0,1]^2$}.
\mg{The problem is governed by the Reynolds-averaged Navier--Stokes equations for incompressible flows, with Spalart--Allmaras turbulence model:
\begin{equation}\label{eq:nssa}
    \begin{cases}
        \displaystyle\frac{\partial \vec{u}}{\partial t} + \nabla \cdot (\vec{u} \otimes \vec{u}) - \nabla \cdot (2(\nu+\nu_t)\nabla^{s}\vec{u}) + \nabla p = \vec{0},                                                                                                                                                                                                            & (\x,t) \in \spatialdomain \times (0,T), \\
        \nabla \cdot \vec{u} = 0,                                                                                                                                                                                                                                                                                                                                                 & (\x,t) \in \spatialdomain \times (0,T), \\
        \displaystyle\frac{\partial \widetilde{\nu}}{\partial t} + \nabla \cdot (\vec{u} \widetilde{\nu}) - \frac{1}{\sigma}\nabla \cdot ((\nu+\widetilde{\nu})\nabla \widetilde{\nu}) - \frac{c_{b2}}{\sigma} \nabla \widetilde{\nu} \cdot \nabla \widetilde{\nu} - c_{b1} \widetilde{S}\widetilde{\nu} + c_{w1} f_{w} \left(\frac{\widetilde{\nu}}{\widetilde{d}}\right)^2 = 0, & (\x,t) \in \spatialdomain \times (0,T), \\
    \end{cases}
\end{equation}
where $\vec{u}(\x,t;\y)$, $p(\x,t;\y)$ and $\widetilde{\nu}(\x,t;\y)$ are the parameter-dependent velocity, pressure and eddy viscosity fields.
Moreover, the turbulent viscosity is defined as
\begin{equation}
    \nu_t = \widetilde{\nu} f_{v1},
\end{equation}
and $\nu$, $\nabla^{s}\vec{u}=(\nabla \vec{u} + (\nabla \vec{u})^\top)/2$ and $\widetilde{d}$ respectively denote the physical viscosity of the fluid, the strain-rate tensor and the distance from the closest wall in the domain.
Finally, the Spalart--Allmaras model is completed with a set of closure functions \cite{Spalart1992}:
\begin{equation}\label{eq:sa}
    \begin{aligned}
        \bm{\Omega}   & :=\displaystyle\frac{\nabla \vec{u} - (\nabla \vec{u})^\top}{2},                                          &
        \widetilde{S} & :=[2\bm{\Omega} : \bm{\Omega}]^{1/2}+\displaystyle\frac{\widetilde{\nu}}{\kappa^2\widetilde{d}^2} f_{v2}, &
        \chi          & := \displaystyle\frac{\widetilde{\nu}}{\nu},                                                                \\
        f_w           & :=g\left[\displaystyle\frac{1+c_{w3}^6}{g^6+c_{w3}^6}\right]^{1/6},                                       &
        f_{v2}        & :=1-\displaystyle\frac{\chi}{1+\chi f_{v1}},                                                              &
        f_{v1}        & :=\displaystyle\frac{\chi^3}{\chi^3 + c_{v1}^3},                                                            \\
        g             & :=r+c_{w2}(r^6-r),                                                                                        &
        r             & :=\frac{\widetilde{\nu}}{\widetilde{S} \kappa^2 \widetilde{d}^2},                                         &
    \end{aligned}
\end{equation}
where $\sigma=2/3$, $c_{b1}:=0.1355$, $c_{b2}=0.622$, $\kappa=0.41$, $c_{w1}=c_{b1}/\kappa^2 + (1+ c_{b2})/\sigma$, $c_{v1}=7.1$, $c_{w2}=0.34$ and $c_{w3}=2$.}

\begin{figure}
    \begin{subfigure}[t]{\textwidth} \centering
        \inputfrom{figures/openfoam/}{schematic.tex}
        \caption{Schematic of boundary conditions}
    \end{subfigure}
    \begin{subfigure}[t]{\textwidth}
        \includegraphics[width=\textwidth,trim=1cm 10cm 1cm 10cm, clip]{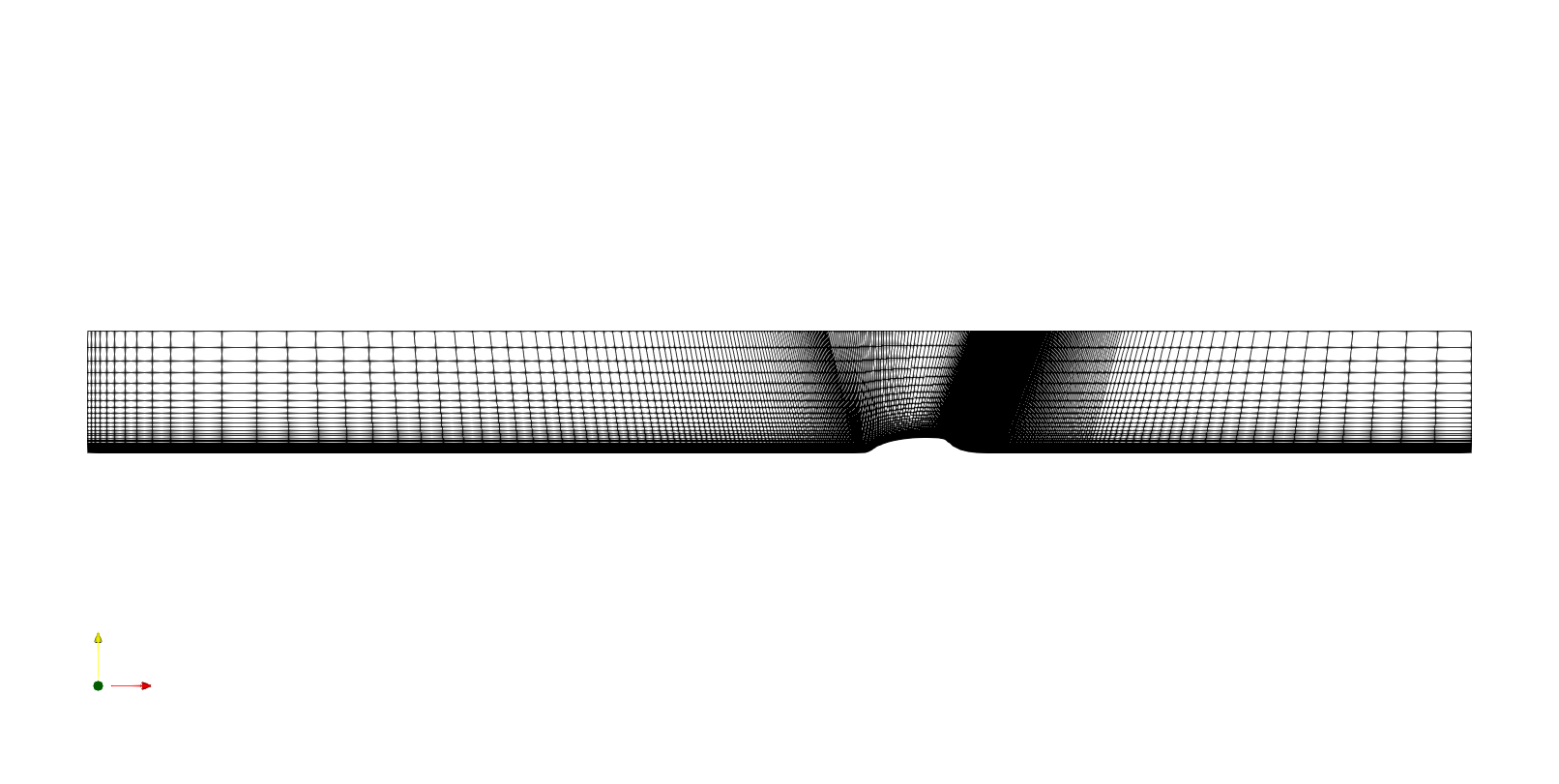}
        \caption{Coarse Mesh}
    \end{subfigure}
    \begin{subfigure}[t]{\textwidth}
        \includegraphics[width=\textwidth,trim=1cm 10cm 1cm 10cm, clip]{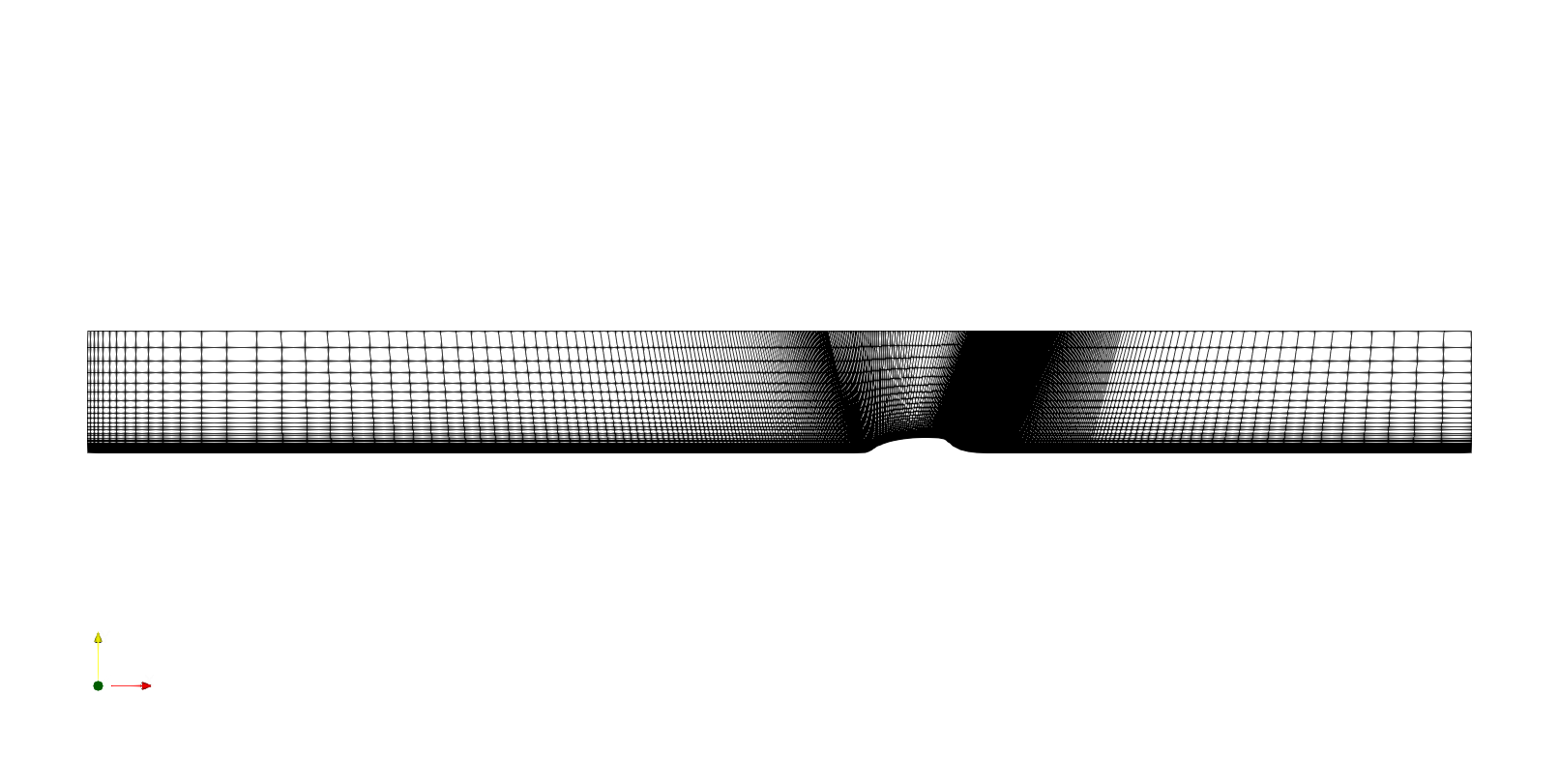}
        \caption{Medium Mesh}
    \end{subfigure}
    \begin{subfigure}[t]{\textwidth}
        \includegraphics[width=\textwidth,trim=1cm 10cm 1cm 10cm, clip]{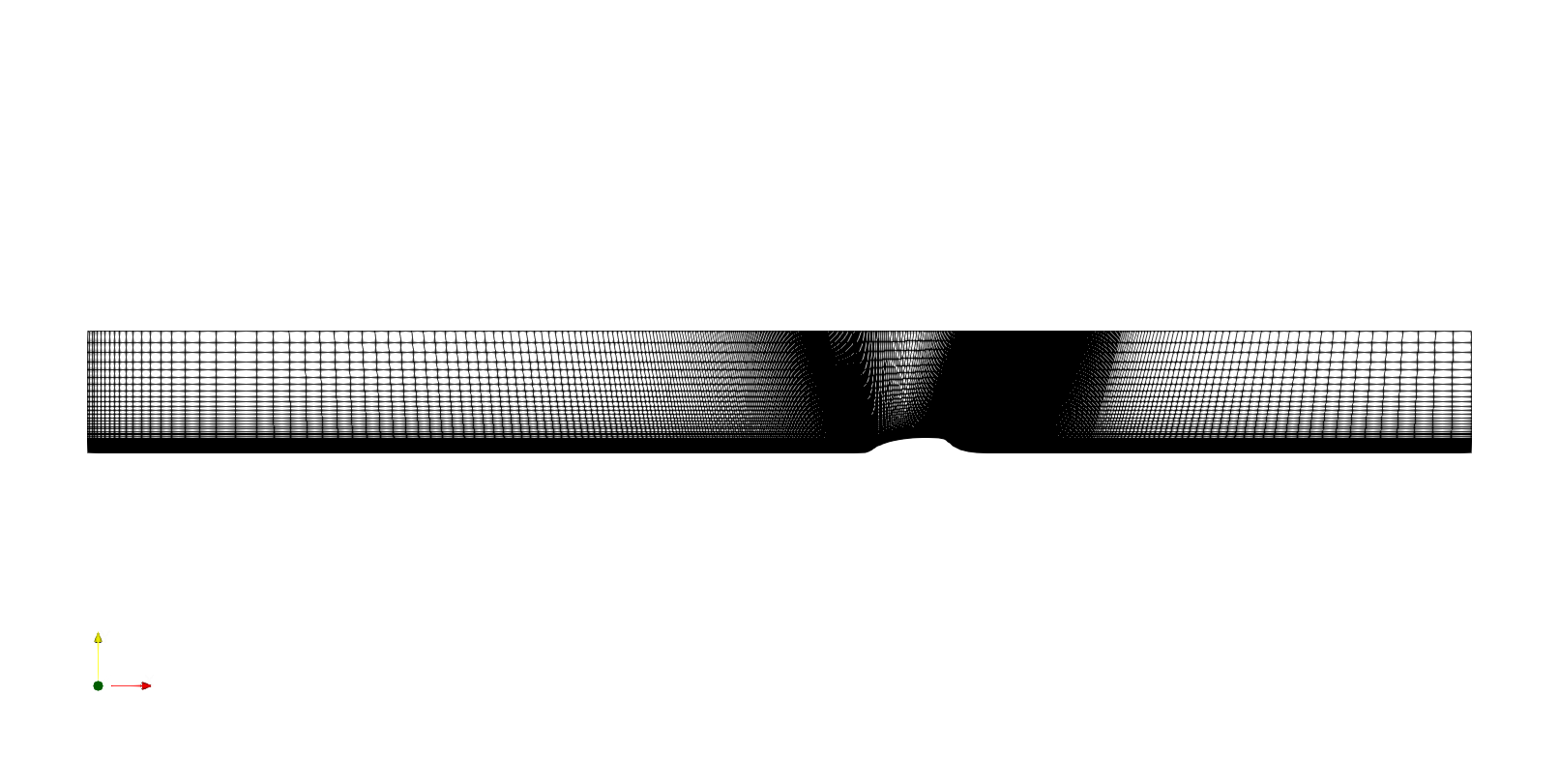}
        \caption{Fine Mesh}
    \end{subfigure}
    \caption{Geometry of spatial domain $\spatialdomain$, \mg{boundary conditions and three examples of} spatial meshes for the NASA 2DWMH test problem described in \Cref{example:2dwmh}.}
    \label{fig:examplemesh}
\end{figure}

The geometry of the spatial domain $\spatialdomain$, \mg{the set of boundary conditions} and three example meshes are shown in \Cref{fig:examplemesh}.
\mg{Let $L=0.42m$ denote the chord length of the hump, with its maximum height set to $H=0.0535m$.
    The computational domain $\spatialdomain$ is a channel of height $L$ extending $6.39L$ upstream and $5L$ downstream of the hump.
    At the inlet $\Gamma_\text{in}$, a parabolic profile is imposed on both velocity and turbulent viscosity, with the maximum velocity magnitude parameterised by $\inflow:\Gamma\to\R$,
    \begin{equation}
        \inflow(\y):=17.3 + (69.2-17.3)y_2,
    \end{equation}
    whereas the eddy viscosity is set to $\widetilde{\nu}=3\nu$, assuming fully developed turbulence.
    A suction jet is located at approximately $65\%$ of the chord length $L$ on a patch named $\Gamma_\text{jet}$.
    On $\Gamma_\text{jet}$, homogeneous Neumann conditions are imposed on the eddy viscosity and a smooth velocity profile is defined, see \cite{Tsiolakis2022}, with maximum outflow magnitude $\jet:\Gamma\to\R$,
    \begin{equation}
        \jet(\y):=10y_1.
    \end{equation}
    Finally, homogeneous Dirichlet conditions are enforced for both $\vec{u}$ and $\widetilde{\nu}$ on $\Gamma_\text{wall}$, homogeneous Neumann conditions are applied on $\Gamma_\text{out}$ and a symmetry condition is applied on $\Gamma_\text{sym}$.

    The kinematic viscosity of the fluid is set to $\nu=1.55274 \cdot 10^{-5} m^2/s$.
    Setting the chord length $L$ as the characteristic length of the problem and the characteristic velocity equal to the maximum inlet velocity, the resulting Reynolds number $\textrm{Re}:=(\inflow L)/\nu$ ranges between $4.70\cdot 10^{5}$ and $1.88\cdot10^{6}$ --- which is clearly in the turbulent regime.
}

\begin{figure}
    \begin{subfigure}{0.49\textwidth}
        \includegraphics[width=\textwidth]{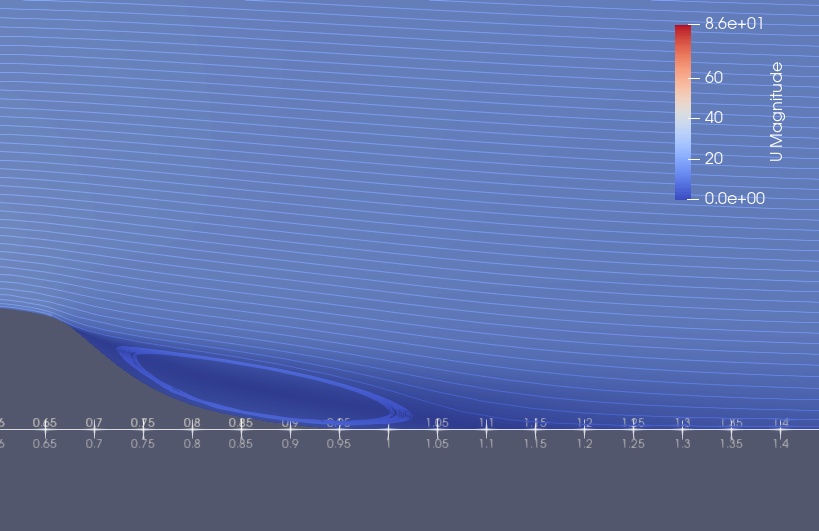}
        \caption{$\y=[1,0]$ i.e.\ $\jet=10$, $\inflow=17.3$}
    \end{subfigure}
    \begin{subfigure}{0.49\textwidth}
        \includegraphics[width=\textwidth]{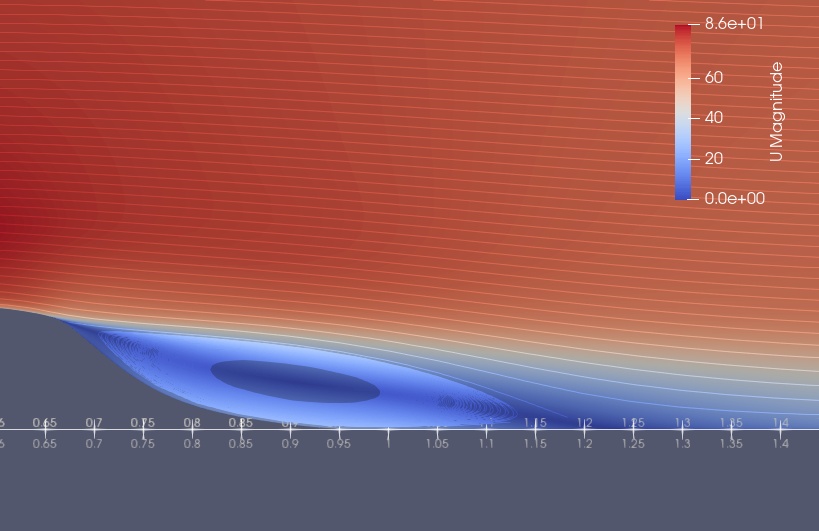}
        \caption{$\y=[0,1]$ i.e.\ $\jet=0$, $\inflow=69.2$}
    \end{subfigure}
    \caption{Velocity streamlines for two parameter realisations $\y=[1,0]$ and $\y=[0,1]$ in the region behind the hump.}
    \label{fig:examplestreamlines}
\end{figure}

\bk{The parametric formulation is motived by a desire to control the flow reattachment point (i.e., the QoI $q:\Gamma\to\R$) by using a suction jet of varying magnitude at the back of the hump.}
\mg{Example streamlines for two flow realisations are shown in \Cref{fig:examplestreamlines}}, \bk{demonstrating the variability in the flow separation and the flow reattachment point.}

\subsubsection{Numerical approximation of the flow fields}
\mg{To compute the reattachment point, \eqref{eq:nssa} is solved using a cell-centered finite volume scheme with pseudo-transient relaxation, and a staggered algorithm splitting the flow equations and turbulence model.
    The \simplefoam{} solver in \openfoam{} \cite{openfoam} provides an implementation of the SIMPLE algorithm \cite{Patankar1972} to simulate such a problem.
    Note that, although the problem under analysis is steady-state, the \simplefoam{} solver employs a pseudo-transient formulation to treat the non-linear terms, thus motivating the presence of the time derivative in \eqref{eq:nssa}.}

\ah{The reattachment point is postprocessed by evaluating the wall shear stress along the lower boundary: the point where the shear stress changes sign indicates reattachment.
    This computation requires accurate velocity fields, especially near the wall}.
\mg{Indeed, it is well known that a value of the \ah{dimensionless wall distance} $y^{+}<1$ is required in order to capture the physical phenomena in the vicinity of the walls.
    In the following subsection, it will be shown that the \plateaumisc{} algorithm can extract relevant information from coarse meshes ($y^{+} \gg 1$), normally unsuitable for turbulent simulations, automatically identifying when the potential of a mesh has been exhausted.}
\bk{These under-resolved meshes are expected to provide quantitatively inaccurate, and possibly noisy, predictions of the quantity of interest}.


\subsubsection{Multi-fidelity MISC approximation}
Three fidelities are defined using the coarse, medium and fine meshes \ah{(see \Cref{fig:examplemesh}), each with distinct horizontal and vertical cell grading, and solver tolerances}.
The vertical \texttt{blockMesh} grading \ah{details are listed in} \Cref{table:blockmeshparams}.
\mg{\begin{table}\centering
        \caption{\openfoam{} \texttt{blockMesh} configuration parameters for the coarse, medium and fine meshes in the 2DWMH problem in \Cref{example:2dwmh}.}
        \label{table:blockmeshparams}
        \begin{tabular}{l r r r}
            Fidelity                                           & Coarse              & Medium               & Fine                \\
            $\miModel$                                         & $1$                 & $2$                  & $3$                 \\
            \hline
            Channel height ($L$)                               & $1$                 & $1$                  & $1$                 \\
            Number of cells in vertical direction              & $42$                & $48$                 & $ 96$               \\
            Height of first cell                               & $4.25\times10^{-4}$ & $1.83\times10^{-4}$  & $1.00\times10^{-5}$ \\
            Height of last cell                                & $1.31\times10^{-1}$ & $ 1.31\times10^{-1}$ & $9.16\times10^{-2}$ \\
            Total expansion ratio (last height / first height) & $308$               & $712$                & $9156$              \\
            Cell-to-cell expansion ratio                       & $1.15$              & $1.15$               & $1.10$              \\
            Range of $y^{+}$ for $\inflow=69.2$, $\jet=0$      & $[6.96,23.78]$      & $[2.41,13.12]$       & $[0.09,0.99]$       \\
        \end{tabular}
    \end{table}}
\mg{The mesh size at the wall and the expansion ratio play a critical role in the accuracy of the solution and, consequently, in the computation of the reattachment point.
    As previously mentioned, a value of the dimensionless height of the first boundary cell, $y^{+}$, below $1$ is required to fully resolve the boundary layer without using wall functions}.
\Cref{fig:yplus} shows that this is only verified for the fine mesh.
\begin{figure}
    \inputfrom{figures/openfoam/}{yplus.tex}
    \caption{Computed $y^+$ values for the coarse, medium and fine meshes in the 2DWMH problem in \Cref{example:2dwmh}. The coordinate $x/L$ is the scaled $x$ coordinate along the bottom boundary where $L$ is the hump length and the origin is the start of the hump.}
    \label{fig:yplus}
\end{figure}
\mg{Hence,} one would expect that only the fine mesh can accurately compute the reattachment point.
\mg{Nonetheless, under-resolved meshes can be employed by \plateaumisc{} to extract some useful information from the problem}.
The solver tolerances are \ah{listed} in \Cref{table:configparams}.
\begin{table}\centering
    \caption{\openfoam{} \simplefoam{} configuration parameters for the coarse, medium and fine meshes in the 2DWMH problem in \Cref{example:2dwmh}.}
    \label{table:configparams}
    \begin{tabular}{l c  c  c}
        Fidelity                                    & Coarse    & Medium    & Fine             \\
        $\miModel$                                  & $1$       & $2$       & $3$              \\
        \hline
        Cells                                       & 18060     & 24384     & 69120            \\
        \hline
        Tolerances (residual)                       &           &           &                  \\
        velocity                                    & $10^{-3}$ & $10^{-4}$ & $10^{-6}$        \\
        turbulence                                  & $10^{-3}$ & $10^{-4}$ & $10^{-6}$        \\
        pressure                                    & $10^{-2}$ & $10^{-3}$ & $5\times10^{-6}$ \\
        \hline
        Tolerances (inner loop solver)              &           &           &                  \\
        absolute                                    & $10^{-6}$ & $10^{-6}$ & $10^{-10}$       \\
        relative                                    & $10^{-2}$ & $10^{-2}$ & $10^{-4}$        \\
        \hline
        Computational cost (sample mean solve time) & 36.9 s    & 183.1 s   & 5062.8 s         \\
    \end{tabular}
\end{table}
A combination of the mesh size and iterative solver tolerances controls the solver error.
In \Cref{fig:hist-openfoam} \ah{approximation error} histograms for the \mg{fidelities} $\miModel=1,2$ are shown, with errors computed relative to the \mg{finest fidelity} $\miModel=3$.
This shows that the solver error decreases as the fidelity $\miModel$ is increased, \ah{although} no clear distribution \ah{emerges}.
\lt{Note also that the histograms are not centered around zero, i.e.\ both fidelities have a non-negligible bias.}

\begin{figure}
    \begin{subfigure}{0.49\textwidth}
        \inputfrom{figures/openfoam/}{histogram_1.tex}
        \caption{Coarse fidelity $\miModel=1$}
    \end{subfigure}
    \begin{subfigure}{0.49\textwidth}
        \inputfrom{figures/openfoam/}{histogram_2.tex}
        \caption{Medium fidelity $\miModel=2$}
    \end{subfigure}
    \caption{Histograms of approximation error $q^\text{ref} - q^{\miModel}$ from $n=1000$ Latin Hypercube samples for the 2DWMH problem using fidelities $q^{\miModel}$ defined in \Cref{table:configparams}}
    \label{fig:hist-openfoam}
\end{figure}

\ah{To quantify} computational cost \ah{for each fidelity, we average the solver runtime over} a sample of parameter realisations\lt{: specifically} samples from a sparse grid of Leja points with the two-step rule from \eqref{eq:rules} and isotropic multi-index set $\{\vMiParam \in\N^2 \st \Vert \vMiParam \Vert_1 \leq 2 + w\}$ with $w=6$.
\ah{All simulations were performed on} a workstation with an Intel Xeon Silver 4215R CPU @ 3.20GHz.
\ah{\Cref{table:configparams} presents} the estimated costs: the \mg{lowest fidelity} is \ah{roughly} five times \ah{faster} than the intermediate fidelity, and 137 times \ah{faster} than the finest fidelity.

\ah{The QoI (the reattachment point) is approximated} using MISC and \plateaumisc{}, \mg{see} \Cref{alg:misc,alg:pmisc} respectively.
A reference approximation is computed with isotropic Smolyak multi-index set \eqref{eq:isotropic-miset} with $w=5$, Clenshaw--Curtis points and the doubling rule, and the fine fidelity $\miModel=3$ defined in \Cref{table:configparams}.
The final response surfaces are shown in \Cref{fig:surfaces2dwmhmisc}.
The \plateaumisc{} approximation is \ah{clearly} more accurate \ah{and less corrupted by solver error than the classical MISC strategy}.
\begin{figure}
    \begin{subfigure}[t]{0.48\textwidth}
        \inputfrom{figures/openfoam/misc_leja/}{surface_noplateau.tex}
        \caption{\textsc{MISC} response surface $\approxFinal{q}(\y)$.}
        \label{fig:error2dwmhmiscsurf}
    \end{subfigure}
    \begin{subfigure}[t]{0.48\textwidth}
        \inputfrom{figures/openfoam/misc_leja/}{surface.tex}
        \caption{\plateaumisc{} response surface $\approxFinal{q}(\y)$.}
        \label{fig:2dwmhmiscsurf}
    \end{subfigure}
    \begin{subfigure}[t]{0.48\textwidth}
        \inputfrom{figures/openfoam/misc_leja/}{errorsurface_noplateau.tex}
        \caption{\textsc{MISC} error surface $\approxRef{q}(\y)-\approxFinal{q}(\y)$.}
        \label{fig:errorsurfopenfoamclassicmisc}
    \end{subfigure}
    \begin{subfigure}[t]{0.48\textwidth}
        \inputfrom{figures/openfoam/misc_leja/}{errorsurface.tex}
        \caption{\plateaumisc{} error surface $\approxRef{q}(\y)-\approxFinal{q}(\y)$.}
        \label{fig:errorsurfopenfoammisc}
    \end{subfigure}
    \caption{Approximated response surfaces $\approxFinal{q}(\y)$ for the 2DWMH problem using \mg{(a) \Cref{alg:misc} and (b) \Cref{alg:pmisc}} with \mg{fidelities} $q^{\miModel}$ defined in \Cref{table:configparams}.
        \Cref{fig:errorsurfopenfoammisc,fig:errorsurfopenfoamclassicmisc} illustrate the \mg{corresponding} approximation errors $\approxRef{q}(\y)-\approxFinal{q}(\y)$.}
    \label{fig:surfaces2dwmhmisc}
\end{figure}

The $L_\rho^2(\Gamma;R)$, and similarly the $H_\rho^1(\Gamma;\R)$, norms of the approximation error are estimated \bk{numerically via} \eqref{eq:estimated-L2} with $N_{MC}=10^4$ Monte Carlo samples.
\ah{These are presented} in \Cref{fig:2dwmherror}.
\begin{figure}
    \begin{subfigure}[t]{0.48\textwidth}
        \inputfrom{figures/openfoam/misc_leja/}{ltwo.tex}
        \caption{Estimated $L_\rho^{2}(\Gamma;\R)$ error}
    \end{subfigure}
    \begin{subfigure}[t]{0.48\textwidth}
        \inputfrom{figures/openfoam/misc_leja/}{hone.tex}
        \caption{Estimated $H_\rho^{1}(\Gamma;\R)$ error}
    \end{subfigure}
    \caption{Estimated $L_{\rho}^2(\Gamma;\R)$ and $H_{\rho}^1(\Gamma;\R)$ approximation error for the \plateaumisc{} approximation $\approxFinal{q}$ for the 2DWMH problem using \Cref{alg:pmisc} with \mg{fidelities} $q^{\miModel}$ defined in \Cref{table:configparams}. Vertical arrows denote the times at which a fidelity is marked as saturated.}
    \label{fig:2dwmherror}
\end{figure}
At low computational budgets, both approximations are poor but are informative of underlying trends and approximate values; note that at such a low cost, \lt{no approximation is attainable} using only the high fidelity solver.
As the computational budget increases, the \plateaumisc{} approximation error decreases at a rate comparable to the reference approximation\lt{, whereas the} MISC approximation error stagnates and does not give a useful approximation, with large errors even at large computational budgets.
This \bk{again} demonstrates that \plateaumisc{} is robust to the solver errors \mg{that are hindering the accuracy of the classical MISC approximation}. 
\Cref{alg:pmisc} is able to balance solver cost and approximation error to automatically increase the fidelity and extract as much information as possible for a given computational budget.

Sampling profiles and spectral envelopes for the \plateaumisc{} approximation are shown in \Cref{fig:2dwmhsamplesenvelopes}.
\begin{figure}
    \begin{subfigure}[t]{0.48\textwidth}
        \inputfrom{figures/openfoam/}{sampleprofile-pmisc.tex}
        \caption{Number of solves $\nColloc{}\rvert_{\miModel}$ for each \mg{fidelity} $\miModel$ at each iteration.}
        \label{fig:2dwmhsamplingprofile-pmisc}
    \end{subfigure}
    \begin{subfigure}[t]{0.48\textwidth}
        \inputfrom{figures/openfoam/}{plateaudetection_env_pmisc.tex}
        \caption{Spectral envelope $\envelope\rvert_{\miModel}$ for each \mg{fidelity} $\miModel$ at \bkrevise{plateau detection}
            \bkrevise{(which happens at different iterations for different fidelities)}.}
        \label{fig:2dwmhenvelope-pmisc}
    \end{subfigure}
    \caption{Number of \mg{fidelity} solves against computational work and snapshot of spectral polynomial coefficients for the \plateaumisc{} approximation $\approxFinal{q}$ of the 2DWMH prolem using \Cref{alg:pmisc} with \mg{fidelities} $q^{\miModel}$ defined in \Cref{table:configparams}.}
    \label{fig:2dwmhsamplesenvelopes}
\end{figure}
\Cref{fig:2dwmhsamplingprofile-pmisc} in combination with \Cref{fig:2dwmherror} illustrates that substantial error reduction only occurs when the number of samples from the high-fidelity ($\miModel=3$) increases.
Spectral envelopes in \Cref{fig:2dwmhenvelope-pmisc} (despite not having theoretical proof that the QoI is analytic with respect to the parameters) exhibit an exponential decay followed by spectral plateaus.
This allows \plateaumisc{} to successfully identify stagnation and \ah{avoid overfitting to the solver noise}, unlike classical MISC.

\mg{It is worth highlighting that \plateaumisc{} is capable of using few samples from under-resolved meshes with $y^{+}=\bigO{10^{1}}$ to extract relevant features of the flow, without ``overfitting'' the solver error.
    This allows us to improve the accuracy of the multi-fidelity surrogate (see \Cref{fig:surfaces2dwmhmisc}), and also significantly reduces the evaluations of all configurations.
    The number of evaluations of fidelity $\miModel=3$, the only fidelity with $y^{+}<1$, is reduced from $41$ evaluations in classic MISC to $27$ evaluations in \plateaumisc{}}.
\ah{The final multi-index sets, shown in \Cref{fig:miset-2dwmh}, also illustrate the tendency to oversample the lower fidelities by the classic MISC algorithm.}
\begin{figure}
    \begin{subfigure}{0.49\textwidth}
        \inputfrom{figures/openfoam/misc_leja/}{miset_np.tex}
        \caption{Multi-index set $\miset$ for final MISC approximation}
    \end{subfigure}
    \begin{subfigure}{0.49\textwidth}
        \inputfrom{figures/openfoam/misc_leja/}{miset.tex}
        \caption{Multi-index set $\miset$ for final \plateaumisc{} approximation}
    \end{subfigure}
    \caption{Multi-index sets $\miset$ for the MISC and \plateaumisc{} approximation $\approxFinal{q}$ of the 2DWMH prolem using \Cref{alg:misc,alg:pmisc} with fidelities $q^{\miModel}$ defined in \Cref{table:configparams}.}
    \label{fig:miset-2dwmh},
\end{figure}

\ah{Finally, \Cref{fig:2dwmh-pdf} compares estimates of the probability density function of the QoI, constructed using $N_{MC}=10^4$ Monte Carlo samples.
    The MISC approximation deviates significantly from the reference, while \plateaumisc{} matches it closely}.
\ah{\Cref{fig:2dwmh-ks2} shows the evolution of the KS2 test statistic which stagnates for the MISC approximation, but decreases consistently for \plateaumisc{} illustrating its superior performance}.
\begin{figure}
    \begin{subfigure}[t]{0.49\textwidth}
        \inputfrom{figures/openfoam/misc_leja/}{pdf.tex}
        \caption{Approximate probability density functions for the final approximations}
        \label{fig:2dwmh-pdf}
    \end{subfigure}
    \begin{subfigure}[t]{0.49\textwidth}
        \inputfrom{figures/openfoam/misc_leja/}{ks.tex}
        \caption{Evolution of two sample Kolmogorov--Smirnov (KS2) test statistic with computational cost}
        \label{fig:2dwmh-ks2}
    \end{subfigure}
    \caption{Final approximate probability density functions and evolution of two sample Kolmogorov--Smirnov test statistic for the \plateaumisc{} approximation $\approxFinal{q}$ of the 2DWMH prolem using \Cref{alg:pmisc} with \mg{fidelities} $q^{\miModel}$ defined in \Cref{table:configparams}.}
\end{figure}

\section{Conclusions} \label{sec:conclusions}
Using the standard adaptive multi-index stochastic collocation algorithm \mg{to compute surrogates} of parametric PDEs can result in an approximation whose quality stagnates since the algorithm can not distinguish between the solver \mg{error} and the true parametric response that must be captured.
\mg{The} \plateaumisc{} algorithm can identify the effect of solver \bk{noise} by inspecting the spectral content of the constructed MISC approximation: a plateau detection algorithm identifies and blocks further refinement in the saturated fidelities, resulting in a multi-fidelity MISC approximation robust to \mg{such intrinsic noise}.

The accuracy of the \plateaumisc{} algorithm is demonstrated through an analytic Gaussian peak example, a parabolic PDE example \mg{with uncertain coefficients} and a \mg{parametric} turbulent Navier--Stokes flow.
In all of these examples the \plateaumisc{} algorithm is able to converge at a rate comparable to the high fidelity reference approximation.
In the analytic Gaussian example and the parametric parabolic PDE \bk{example}, a significant reduction in computational cost is observed when compared to a single high fidelity approximation.
In the \mg{parametric} turbulent flow example, \mg{\plateaumisc{} numerically demonstrates the ability to automatically select the appropriate fidelities required to construct a robust surrogate model for quantities of engineering interest, whilst extracting relevant information even from under-resolved meshes not suitable for reliable single fidelity computations.}
In other words, the advantage of multi-fidelity approximation over single-fidelity is dependent on the pool of available fidelities and
\plateaumisc{} is able to exploit whatever advantage is there due to its robustness to noise, whilst classic MISC may result in suboptimal approximations overfitting intrinsic solver errors.


\MGnew{
Finally, we remark that this work has proposed a robust variant of MISC, outperforming the standard one in the presence of noisy evaluations of quantities of interests.
Whilst we can point out that the reduced number of hyperparameters to be tuned represents an appealing feature of \plateaumisc{} also with respect to other approaches (e.g., based on least-squares techniques), a thorough benchmarking is required to globally assess advantages and disadvantages of the different noise-robust multi-fidelity strategies in the literature.
In this context, a benchmarking campaign in the spirit of \cite{Seelinger2025,Jakeman2025} is expected to significantly contribute to a deeper understanding of the performance of noise-robust multi-fidelity methods, including \plateaumisc{}.
For this purpose, a crucial point is represented by the selection of appropriate test problems accounting for the different aspects possibly influencing each method, including, e.g., choice of samples, choice of polynomial spaces, and algorithms to allocate over fidelities, just to name a few.}

\section{Acknowledgments}
B.\ Kent and L.\ Tamellini have been partially supported by the project 202222PACR “Numerical approximation of uncertainty quantiﬁcation problems for PDEs by multi-ﬁdelity methods (UQ-FLY)”, funded by European Union--NextGenerationEU.
L.\ Tamellini has been partially supported by the project PRIN PNRR “Uncertainty Quantification of coupled models for water flow and contaminant transport” (P2022LXLYY), financed by the European Union--NextGeneration EU.
L.\ Tamellini is member of the Gruppo Nazionale Calcolo Scientiﬁco-Istituto Nazionale di Alta Matematica (GNCS-INdAM).

\mgrevise{M. Giacomini and A. Huerta acknowledge the support of the Spanish Ministry of Science, Innovation and Universities and Spanish State Research Agency MICIU/AEI/10.13039/501100011033 (Grant No. PID2023-149979OB-I00) and of the Generalitat de Catalunya (Grant No. 2021- SGR-01049). M. Giacomini is Fellow of the Serra Húnter Programme of the Generalitat de Catalunya.}

\bibliographystyle{elsarticle-num}

\end{document}